\documentclass[a4paper,fleqn]{cas-sc}
\usepackage[numbers]{natbib}
\usepackage[utf8]{inputenc}
\usepackage[T1]{fontenc}
\usepackage{amsmath,amssymb,mathtools}
\usepackage{booktabs}
\usepackage{graphicx}
\usepackage{microtype}
\usepackage{url}
\usepackage{placeins}
\usepackage{float}
\restylefloat{figure}
\hypersetup{hypertexnames=false}
\numberwithin{equation}{section}
\graphicspath{{./}}

\begin{document}
\let\WriteBookmarks\relax
\shorttitle{Conservative Local Neural Operator (CLNO)}
\shortauthors{B. Wang, R. Li and Q. Zou}
\title[mode=title]{CLNO: Conservative Local Neural Operator for Compressible Flow Prediction}
\author[1]{Buxuan Wang}
\ead{wangbx29@mail2.sysu.edu.cn}
\credit{Conceptualization, Methodology, Software, Validation, Formal analysis, Investigation, Visualization, Writing -- original draft}
\author[1]{Rongzu Li}
\ead{lirz6@mail2.sysu.edu.cn}
\author[1]{Qingsong Zou}
\cormark[1]
\ead{mcszqs@mail.sysu.edu.cn}
\credit{Supervision, Project administration, Funding acquisition, Writing -- review and editing}
\affiliation[1]{organization={School of Computer Science and Engineering, Sun Yat-sen University},
                city={Guangzhou}, country={China}}
\cortext[1]{Corresponding author.}

\begin{abstract}
We propose a conservative local neural operator (CLNO) for accurate and generalizable autoregressive prediction of compressible flows. Building upon the locality and translation invariance of the local neural operator (LNO), CLNO introduces a conservative evolution framework based on local subdomains of fixed physical size, separating local mean-state evolution from fine-scale reconstruction. The flow is encoded into subdomain-integrated conservative states and fine-scale features, while a shared local network predicts interface fluxes for finite-volume updates that enforce discrete conservation. A convolutional reconstruction network recovers fine-scale structures using zero-mean basis functions, preserving the updated local means and providing information for subsequent evolution. To accommodate different physical parameters, Reynolds--Mach number conditioning modulates the encoding, flux-prediction, and reconstruction networks, enabling a single model to learn parameter-dependent dynamics and generalize to unseen parameter combinations without retraining. CLNO supports prediction across spatial resolutions, computational domain sizes, and boundary conditions. Numerical experiments demonstrate improved generalization and competitive or higher accuracy compared with LNO, including order-of-magnitude error reductions in representative cases, at approximately three to five times the inference cost on periodic benchmarks. These results demonstrate accurate, generalizable, and physically consistent flow prediction through conservative local evolution and fine-scale reconstruction.
\end{abstract}

\maketitle

\section{Introduction}
\label{sec:introduction}

Accurate prediction of unsteady compressible flows is an important
task in computational fluid dynamics. Such flows involve nonlinear
convection, pressure-wave propagation, viscous dissipation, and
interactions across multiple spatial scales. Conventional numerical
methods, including finite-volume and finite-difference schemes,
solve the governing equations through spatial discretization and
successive time integration \cite{leveque2002finite}. However,
high-resolution simulations can be computationally expensive,
particularly for long-time prediction and repeated calculations
under different initial conditions, computational domains, and
boundary conditions. This motivates the development of efficient
and generalizable data-driven methods that retain essential
physical properties.

Neural operators provide an alternative approach by learning
mappings between functions from simulation data. DeepONet
represents nonlinear operators using branch and trunk networks
\cite{lu2021deeponet}, while the Fourier neural operator (FNO)
parameterizes integral operators in Fourier space
\cite{li2021fno}. Subsequent studies have extended neural
operator learning to more flexible spatial discretizations
and geometries. Message-passing neural PDE solvers capture
local interactions through graph-based representations
\cite{brandstetter2022message}. The geometry-informed neural
operator (GINO) maps information between irregular physical
discretizations and a shared computational representation
\cite{li2023gino}, while Transolver learns PDE solutions
on general geometries through physics-aware representations
\cite{wu2024transolver}. These methods improve the flexibility
of neural PDE solvers, but their operator representations
do not inherently enforce finite-volume conservation during
autoregressive time evolution.

The local neural operator (LNO) introduces locality and
translation invariance into transient PDE prediction, allowing
learned local evolution rules to be reused across computational
domains with appropriate boundary treatments
\cite{li2023lno}. Further studies have shown that its
prediction accuracy depends strongly on the relationship
between the network receptive range and the physical
propagation range \cite{ye2024locality}. In common
grid-based implementations, the receptive range is
defined by a fixed number of grid points, so its physical
extent changes with spatial resolution. In addition,
truncating the local Legendre representation may limit
the recovery of fine-scale flow structures, particularly
when the receptive range increases. More importantly,
the original LNO directly predicts the next flow state
without explicitly enforcing conservative exchanges
between neighboring regions. Consequently, its local
time-marching structure does not guarantee the discrete
conservation of mass, momentum, and total energy.

Incorporating neural networks into numerical discretizations
provides another approach to physically consistent PDE
prediction. Bar-Sinai et al. developed data-driven spatial
discretizations for coarse-resolution simulations
\cite{barsinai2019discretizations}, while Kochkov et al.
combined learned numerical components with conventional
fluid solvers \cite{kochkov2021mlcfd}. The conservative
form network (CFN) learns numerical flux functions within
conservative updates for hyperbolic conservation laws
\cite{chen2024cfn}. FluxGNN incorporates finite-volume-inspired
conservation and symmetry properties into graph neural PDE
solvers \cite{horie2024fluxgnn}, while entropy-stable
conservative flux-form neural networks introduce additional
constraints on learned conservative dynamics
\cite{liu2024escfn}. FluxNet further learns constrained
cumulative transport between neighboring cells and handles
nonperiodic boundaries using ghost cells
\cite{lan2026fluxnet}. These approaches incorporate physical
constraints through learned discretizations, conservative
fluxes, or graph-based updates. In contrast, our approach
focuses on performing conservative time evolution on local
subdomains of fixed physical size, independently of the
input-grid resolution, while reconstructing fine-scale
flow structures within these subdomains without modifying
the updated local conservative states.

The central objective of this work is to develop a
local neural operator that supports transfer across
spatial discretizations, computational domains, and
boundary conditions while preserving finite-volume
conservation during autoregressive time evolution.
Achieving this requires a local representation defined
in physical space, conservative state updates through
shared interface fluxes, and consistent mappings
between physical-grid fields and local-subdomain states.
An additional challenge arises from local-subdomain
evolution: local mean states alone cannot fully describe
the internal gradients, waves, and vortical structures
that influence flow dynamics. Therefore, fine-scale
reconstruction must recover the unresolved spatial
information without altering the conservative states
determined by the finite-volume update.

Physical-parameter variation presents a further challenge.
For the nondimensional compressible Navier--Stokes equations,
the Reynolds and Mach numbers affect viscous transport and
compressibility, respectively. Even with identical initial
conditions, computational domains, and boundary conditions,
different Reynolds--Mach number combinations generally
produce different flow evolution. Training a separate
neural operator for each combination is therefore
inefficient. A unified model can instead learn
parameter-dependent evolution rules from multiple
parameter combinations and use the prescribed
physical parameters to adapt its local representations,
providing predictive generalization to combinations
not encountered during training.

To address these challenges, we propose a conservative
local neural operator (CLNO) consisting of three
components. First, a local encoder represents the input
flow using subdomain-integrated conservative quantities
and learned fine-scale features on local subdomains of
fixed physical size. Second, a shared convolutional
network predicts one flux at each interface, and a
finite-volume scheme updates the integrated quantities
using equal-and-opposite flux contributions from
neighboring subdomains. Third, a local convolutional
reconstruction network predicts coefficients of
discrete zero-mean Legendre basis functions to recover
fine-scale flow fields within each local subdomain while
preserving the updated local mean states. The conservative
evolution and fine-scale reconstruction are thus both
performed on the local subdomains, whereas the predicted
flow can be evaluated on the requested physical grid.
This formulation permits the same learned interior
evolution rule to be applied across different
physical-grid resolutions and domain sizes, with
additional numerical treatments imposed for
nonperiodic boundaries and obstacles. The shared
interface-flux construction guarantees discrete
conservation of internal transfers, and yields
global conservation on closed periodic domains.
In addition, a physical-parameter conditioning network
encodes the Reynolds and Mach numbers and modulates
the intermediate features of all three components,
enabling shared model parameters across different
flow regimes and predictive transfer to unseen
parameter combinations.

The main contributions of this work are summarized
as follows:

\begin{enumerate}

\item We propose a conservative local neural operator
based on local subdomains of fixed physical size.
By separating subdomain-integrated conservative
quantities from fine-scale features and updating
them through shared interface fluxes, the proposed
method combines local operator learning with
finite-volume conservation and supports transfer
across physical-grid resolutions.

\item We develop a mean-preserving fine-scale
reconstruction method based on discrete zero-mean
Legendre basis functions. The reconstruction
recovers fine-scale flow structures within each
local subdomain while maintaining consistency
with the updated local conservative states,
enabling high-resolution predictions without
violating the prescribed local mean states.

\item We evaluate CLNO across different spatial
resolutions, computational domain sizes, boundary
configurations, and long autoregressive rollouts.
The results demonstrate improved prediction
accuracy and generalization compared with LNO
on the evaluated benchmarks, including
order-of-magnitude error reductions in
representative cases, with moderate additional
inference cost for periodic flow predictions.

\item We introduce physical-parameter modulation
based on the Reynolds and Mach numbers.
A shared conditioning network modulates the local
encoding, interface-flux prediction, and reconstruction
networks, enabling a single model to represent
different physical-parameter combinations and
achieve a degree of predictive generalization
to unseen combinations without retraining.

\end{enumerate}

\section{Problem formulation}
\label{sec:problem}

We consider the nondimensional compressible Navier--Stokes
equations on a computational domain
$\Omega\subset\mathbb R^d$, where $d\in\{1,2,3\}$.
The governing equations are written in conservative
form as
\begin{equation}
\frac{\partial U}{\partial t}
+\nabla\cdot F^c(U;\boldsymbol\mu)
=\nabla\cdot F^v(U,\nabla U;\boldsymbol\mu),
\label{eq:governing}
\end{equation}
where $U(x,t)=(\rho,\rho u_1,\ldots,\rho u_d,E)^{\mathsf T}
\in\mathbb R^{d+2}$ denotes the conservative state,
$F^c$ and $F^v$ are the convective and
viscous fluxes, respectively, and
$\boldsymbol\mu=(\mathrm{Re},\mathrm{Ma})$
contains the Reynolds and Mach numbers.
Their dependence is determined by the adopted
nondimensionalization and constitutive relations.
The initial and boundary conditions are
$U(x,0)=U_0(x)$ and $\mathcal B(U)=0$.

Let $U^n=U(\cdot,t_n)$ and
$t_{n+1}=t_n+\Delta t$. The exact
time-evolution operator $\mathcal G$ satisfies
\begin{equation}
\begin{lgathered}
U^{n+1}
=\mathcal G(U^n;\Omega,\mathcal B,\boldsymbol\mu)
=U^n-\int_{t_n}^{t_{n+1}}
\nabla\cdot\left[
F^c(U;\boldsymbol\mu)
-F^v(U,\nabla U;\boldsymbol\mu)
\right]\,\mathrm dt,
\end{lgathered}
\label{eq:exact-evolution}
\end{equation}
where $U$ is the solution on $\Omega$,
subject to $U(\cdot,t_n)=U^n$ and
$\mathcal B(U)=0$ on $\partial\Omega$.
For fixed initial and boundary conditions
and computational domain, different values
of $\boldsymbol\mu$ generally yield
different time-evolution operators.

Local conservation is a fundamental property
of the governing equations. Consider a
partition of $\Omega$ into nonoverlapping
local subdomains of prescribed physical scale:
\begin{equation}
\mathcal K_\Omega=\{K_s\}_{s=1}^{S},
\qquad
\overline{\Omega}
=\bigcup_{s=1}^{S}\overline{K_s},
\qquad
|K_s\cap K_j|=0\quad(s\ne j).
\label{eq:subdomain-partition}
\end{equation}
Integrating the governing equations over
each $K_s$ and applying the divergence theorem
gives the finite-volume integral form
\begin{equation}
\begin{lgathered}
\int_{K_s}
\left[
\mathcal G(U^n;\Omega,\mathcal B,\boldsymbol\mu)
-U^n
\right]\,\mathrm dx=
-\int_{t_n}^{t_{n+1}}\int_{\partial K_s}
\left[
F^c(U;\boldsymbol\mu)
-F^v(U,\nabla U;\boldsymbol\mu)
\right]n_s\,\mathrm dS\,\mathrm dt,
\end{lgathered}
\label{eq:exact-conservation}
\end{equation}
where $n_s$ is the outward unit normal
to $\partial K_s$.

We seek a unified neural operator
$\mathcal F_{\Omega,\mathcal B}
[\mathcal G_\theta(\boldsymbol\mu)]$
that approximates the corresponding PDE
time evolution:
\begin{equation}
\widehat U^{n+1}
=\mathcal F_{\Omega,\mathcal B}
[\mathcal G_\theta(\boldsymbol\mu)](U^n)
\approx
\mathcal G(U^n;\Omega,\mathcal B,\boldsymbol\mu)
=U^{n+1},
\label{eq:operator-prediction}
\end{equation}
where $\mathcal G_\theta$ denotes a shared
parameterized evolution mapping with
trainable parameters $\theta$.
The physical parameters $\boldsymbol\mu$
modulate its local networks. At each time step,
$\mathcal F_{\Omega,\mathcal B}$ specifies
the local subdomain partition associated
with $\Omega$ and incorporates the boundary
conditions $\mathcal B$ into
$\mathcal G_\theta$. The trainable parameters
are shared across different initial conditions,
spatial resolutions, computational domains,
boundary conditions, and Reynolds--Mach
number combinations. The resulting model
supports direct inference at previously
unseen parameter combinations and a degree
of predictive generalization across
physical-parameter settings.

The predicted flow is further required
to satisfy the finite-volume conservation
law on the same partition:
\begin{equation}
\int_{K_s}\widehat U^{n+1}(x)\,\mathrm dx
=\int_{K_s}U^n(x)\,\mathrm dx
-\Delta t\sum_{\gamma\in\Gamma_s}
\varepsilon_{s,\gamma}A_\gamma
\widehat F_\gamma^n,
\label{eq:operator-conservation}
\end{equation}
where $\Gamma_s$ is the set of faces of $K_s$,
$A_\gamma$ is the face measure,
$\varepsilon_{s,\gamma}\in\{-1,1\}$
specifies the face orientation, and
$\widehat F_\gamma^n$ is the time-averaged
numerical flux. Each internal face has
one shared flux with opposite orientations
in adjacent subdomains. The operator
architecture and implementation are
presented in Section~\ref{sec:method}.

\section{Method}
\label{sec:method}

\subsection{Overall framework}
\label{sec:overall-framework}

\begin{figure}[H]
\centering
\includegraphics[width=\linewidth]{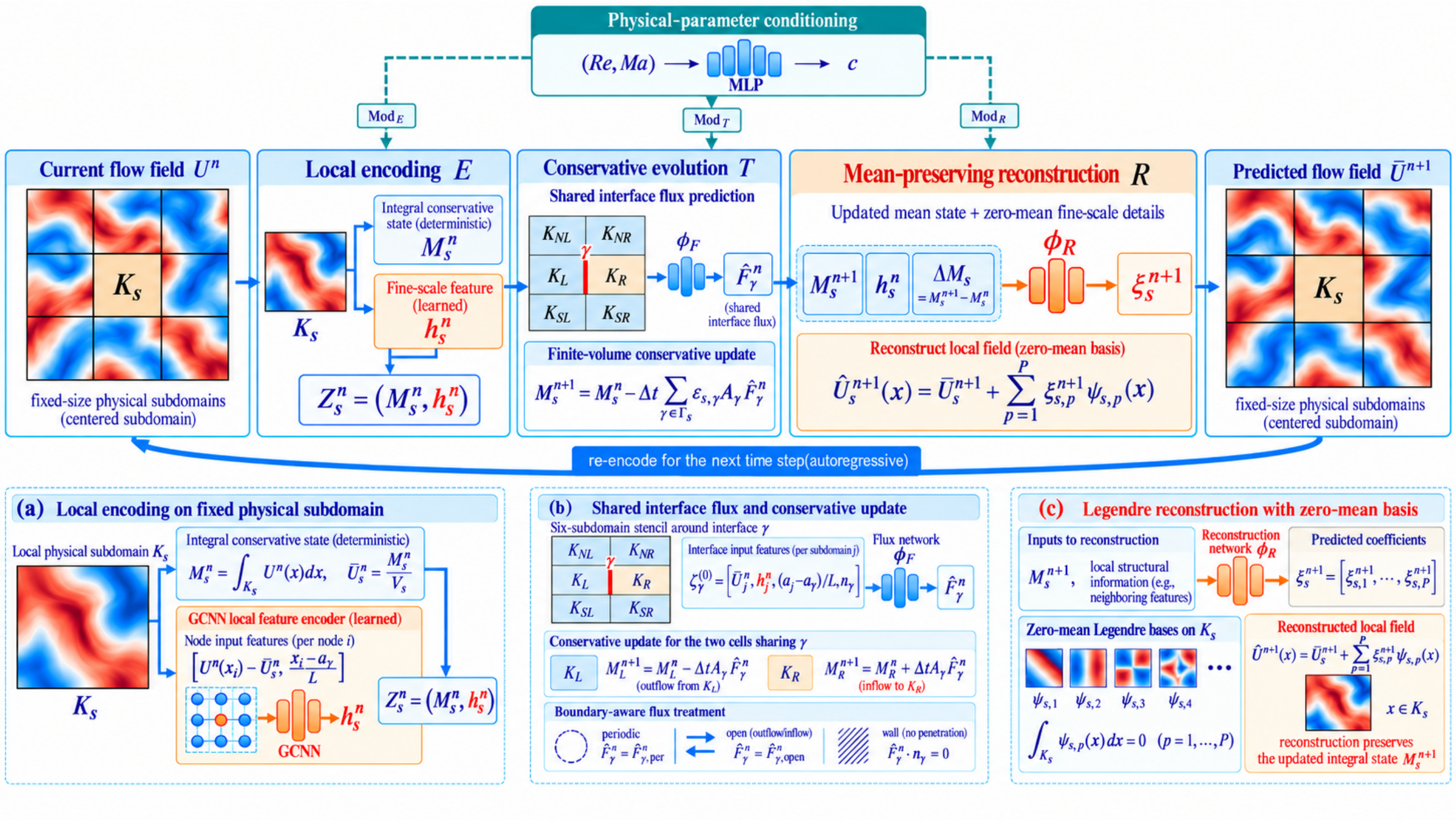}
\caption{Overall framework of CLNO. The current flow field is
encoded on local subdomains of fixed physical size and advanced
through shared interface fluxes. Reconstruction adds zero-mean
fine-scale details to the updated local means, and the predicted
field is re-encoded for the next time step. The Reynolds and Mach
numbers condition all three modules. Panels (a)--(c) detail local
encoding, conservative interface-flux updates, and mean-preserving
Legendre reconstruction, respectively.}
\label{fig:clno-framework}
\end{figure}

As illustrated in Figure~\ref{fig:clno-framework}, the proposed
conservative local neural operator (CLNO) consists of three
modules: local encoding, flux prediction with conservative
evolution, and mean-preserving local reconstruction, denoted by
$\mathcal E$, $\mathcal T$, and $\mathcal R$, respectively.
A physical-parameter conditioning mechanism modulates the three
modules, and the reconstructed flow is re-encoded for the next
autoregressive time step. All trainable parameters are
collectively denoted by $\theta$.

Given the Reynolds number $\mathrm{Re}$ and
Mach number $\mathrm{Ma}$, the physical-parameter
vector $\boldsymbol\mu=(\mathrm{Re},\mathrm{Ma})$
is normalized to obtain the dimensionless conditioning
vector $\mathbf c$.
The Reynolds and Mach numbers are prescribed equation
parameters, independent of the initial condition and
computational-domain size.

The conditioning vector modulates intermediate
features in all three networks through feature-wise
affine transformations. For $X\in\{E,F,R\}$,
the modulation of a selected feature tensor is defined by
\begin{equation}
\begin{lgathered}
\operatorname{Mod}_X(z;\mathbf c)
=\boldsymbol\gamma_X(\mathbf c)\odot z
+\boldsymbol\beta_X(\mathbf c)
\end{lgathered}
\label{eq:parameter-modulation}
\end{equation}
where $X$ identifies the encoding, flux-prediction,
or reconstruction network, and $z$ is an intermediate
feature vector. The symbol $\odot$ denotes
elementwise multiplication, and
$\boldsymbol\gamma_X$ and
$\boldsymbol\beta_X$ are learned mappings
that generate channel-wise scaling and shifting
coefficients of the same dimension as $z$.
The three networks receive the same $\mathbf c$ but use
separate multilayer perceptrons for their modulation mappings.
Because $\mathbf c=\mathbf c(\boldsymbol\mu)$ is computed
entirely from the prescribed physical parameters, we write
the dependence of each module on $\boldsymbol\mu$ alone;
$\mathbf c$ is used only within the modulation mechanism.

The shared parameterized evolution mapping
$\mathcal G_\theta(\boldsymbol\mu)$ is expressed
as the composition
\begin{equation}
\begin{lgathered}
\mathcal G_\theta(\boldsymbol\mu)=\mathcal R[\boldsymbol\mu]\circ\mathcal T[\boldsymbol\mu]\circ\mathcal E[\boldsymbol\mu]
\end{lgathered}
\label{eq:parameterized-evolution}
\end{equation}
where $\circ$ denotes operator composition.

For a computational domain $\Omega$ and
boundary conditions $\mathcal B$, the mapping
$\mathcal F_{\Omega,\mathcal B}$ specifies
the local subdomain partition
$\mathcal K_\Omega=\{K_s\}_{s=1}^{S}$
introduced in Section~\ref{sec:problem}
and incorporates the corresponding boundary
treatment:
\begin{equation}
\begin{lgathered}
\mathcal F_{\Omega,\mathcal B}[\mathcal G_\theta(\boldsymbol\mu)]=\mathcal R[\mathcal K_\Omega,\boldsymbol\mu]\circ\mathcal T[\mathcal K_\Omega,\Delta t,\mathcal B,\boldsymbol\mu]\circ\mathcal E[\mathcal K_\Omega,\boldsymbol\mu]
\end{lgathered}
\label{eq:complete-operator}
\end{equation}
where $\Delta t$ is the time-step size.
The encoding module maps the current flow field
to local states, the evolution module updates
their integral conservative quantities,
and the reconstruction module recovers
the predicted flow field. The resulting
one-step prediction is
\begin{equation}
\begin{lgathered}
\widehat U^{n+1}=\mathcal F_{\Omega,\mathcal B}[\mathcal G_\theta(\boldsymbol\mu)](U^n)
\end{lgathered}
\label{eq:complete-prediction}
\end{equation}
The physical parameters modulate the learned
networks, whereas the computational domain
and boundary conditions determine the
geometric and boundary realization.
All trainable parameters remain shared.

\subsection{Local encoding module}
\label{sec:encoding}

For each local subdomain $K_s$, let
$V_s=|K_s|$ and $a_s$ denote its volume
and center, respectively. The local state
$Z_s^n=(M_s^n,h_s^n)$ consists of the
integral conservative state $M_s^n$
and fine-scale feature vector $h_s^n$:
\begin{equation}
\begin{lgathered}
M_s^n=\int_{K_s}U^n(x)\,\mathrm dx,\qquad \bar U_s^n=\frac{M_s^n}{V_s}
\end{lgathered}
\label{eq:local-state}
\end{equation}
where $\bar U_s^n$ denotes the local
mean conservative state.

A graph convolutional neural network (GCNN)
extracts fine-scale features from input
grid points in a fixed physical
neighborhood of $K_s$. Let
$\mathcal N_E(s)$ denote the indices
of these grid points, with coordinates
$x_i$. Their initial node features are
\begin{equation}
\begin{lgathered}
\eta_{s,i}^{(0)}=\left[U^n(x_i)-\bar U_s^n,\frac{x_i-a_s}{L}\right],\qquad i\in\mathcal N_E(s)
\end{lgathered}
\label{eq:encoder-input}
\end{equation}
where $\eta_{s,i}^{(0)}$ is the initial
feature of node $i$ associated with
subdomain $K_s$, and $L$ is a prescribed
reference length. Graph edges connect
grid points within a fixed physical distance.

At layer $q$, node features are updated
through message passing followed by
physical-parameter modulation:
\begin{equation}
\begin{lgathered}
\eta_{s,i}^{(q+1)}=\operatorname{Mod}_E^{(q)}\left(\operatorname{ReLU}\left[W_E^{(q)}\eta_{s,i}^{(q)}+\frac{1}{|\operatorname{Adj}_s(i)|}\sum_{j\in\operatorname{Adj}_s(i)}\phi_E^{(q)}\left(\eta_{s,i}^{(q)},\eta_{s,j}^{(q)},\frac{x_j-x_i}{L}\right)\right];\mathbf c\right)
\end{lgathered}
\label{eq:graph-convolution}
\end{equation}
where $\eta_{s,i}^{(q)}$ is the feature
at layer $q$, $W_E^{(q)}$ is a trainable
linear transformation,
$\operatorname{Adj}_s(i)$ is the graph
neighborhood of node $i$, including itself,
and $\phi_E^{(q)}$ is a message-generating
multilayer perceptron. The activation
function $\operatorname{ReLU}$ is applied
before parameter modulation.

After $N_E$ graph convolution layers,
weighted pooling produces
\begin{equation}
\begin{lgathered}
h_s^n=\frac{\displaystyle\sum_{i\in\mathcal N_E(s)}\varpi_i\eta_{s,i}^{(N_E)}}{\displaystyle\sum_{i\in\mathcal N_E(s)}\varpi_i}
\end{lgathered}
\label{eq:encoder-pooling}
\end{equation}
where $N_E$ is the number of graph
convolution layers and $\varpi_i$
is the quadrature weight associated
with grid point $x_i$. The encoding
module yields
\begin{equation}
\begin{lgathered}
Z^n=\mathcal E[\mathcal K_\Omega,\boldsymbol\mu](U^n)=(Z_1^n,\ldots,Z_S^n)
\end{lgathered}
\label{eq:encoding-module}
\end{equation}
where $Z^n$ collects the local
states of all subdomains.

\subsection{Flux prediction module}
\label{sec:evolution}

A convolutional neural network (CNN)
predicts the numerical flux across
each internal interface. Let $\gamma$
denote an interface with center
$a_\gamma$ and reference unit normal
$n_\gamma$. Its input features are
\begin{equation}
\begin{lgathered}
\zeta_{\gamma,j}^{(0)}=\left[\bar U_j^n,h_j^n,\frac{a_j-a_\gamma}{L},n_\gamma\right],\qquad j\in\mathcal N_F(\gamma)
\end{lgathered}
\label{eq:flux-input}
\end{equation}
where $\mathcal N_F(\gamma)$ is an
ordered stencil of neighboring subdomains,
$j$ indexes a subdomain in the stencil,
and $\zeta_{\gamma,j}^{(0)}$ denotes
its initial flux-prediction feature.

At layer $q$, the convolutional features
are updated by
\begin{equation}
\begin{lgathered}
\zeta_{\gamma,j}^{(q+1)}=\operatorname{Mod}_F^{(q)}\left(\operatorname{ReLU}\left[\sum_{\delta\in\Lambda_F}W_{F,\delta}^{(q)}\widetilde\zeta_{\gamma,j+\delta}^{(q)}+b_F^{(q)}\right];\mathbf c\right)
\end{lgathered}
\label{eq:flux-convolution}
\end{equation}
where $\Lambda_F$ is the set of
convolution-kernel offsets,
$\delta$ denotes an offset,
$W_{F,\delta}^{(q)}$ and $b_F^{(q)}$
are the trainable convolution weights
and bias, and $\widetilde\zeta$
denotes the padded stencil features.

After $N_F$ convolution layers,
average pooling and linear projection
give the predicted interface flux:
\begin{equation}
\begin{lgathered}
\widehat F_{\gamma,\mathrm{learned}}^n=W_{F,\mathrm{out}}\left(\frac{1}{|\mathcal N_F(\gamma)|}\sum_{j\in\mathcal N_F(\gamma)}\zeta_{\gamma,j}^{(N_F)}\right)+b_{F,\mathrm{out}}
\end{lgathered}
\label{eq:learned-flux}
\end{equation}
where $N_F$ is the number of
flux-prediction convolution layers,
$W_{F,\mathrm{out}}$ and
$b_{F,\mathrm{out}}$ are the output
projection weights and bias, and
$\widehat F_{\gamma,\mathrm{learned}}^n
\in\mathbb R^{d+2}$ denotes the
predicted normal numerical flux,
interpreted as an effective time-averaged
flux over the current time step.

The local integral conservative states
are advanced by the finite-volume scheme
\begin{equation}
\begin{lgathered}
M_s^{n+1}=M_s^n-\Delta t\sum_{\gamma\in\Gamma_s}\varepsilon_{s,\gamma}A_\gamma\widehat F_\gamma^n
\end{lgathered}
\label{eq:fv-update}
\end{equation}
where $\Gamma_s$ is the set of interfaces
of $K_s$, $A_\gamma$ is the interface
measure, and
$\varepsilon_{s,\gamma}\in\{-1,1\}$
specifies the interface orientation
relative to $n_\gamma$.
The numerical flux $\widehat F_\gamma^n$
is selected according to the interface
treatment in Section~\ref{sec:boundary}.

Each internal interface has one shared
flux, which enters the two adjacent
subdomain updates with opposite signs.
Consequently, internal flux contributions
cancel in the discrete global balance.

Let
$\mathbf M^{n+1}
=(M_1^{n+1},\ldots,M_S^{n+1})$
collect the updated integral conservative
states. The evolution module retains
the encoded states and produces
\begin{equation}
\begin{lgathered}
\mathcal T[\mathcal K_\Omega,\Delta t,\mathcal B,\boldsymbol\mu](Z^n)=(Z^n,\mathbf M^{n+1})
\end{lgathered}
\label{eq:evolution-module}
\end{equation}

\subsection{Local reconstruction module}
\label{sec:reconstruction}

The updated local mean state is
$\bar U_s^{n+1}=M_s^{n+1}/V_s$.
A CNN predicts the coefficients of
the local Legendre expansion.
For each subdomain $K_s$, its input
features are
\begin{equation}
\begin{lgathered}
\chi_{s,j}^{(0)}=\left[h_j^n,\bar U_j^n,\bar U_j^{n+1},\frac{a_j-a_s}{L}\right],\qquad j\in\mathcal N_R(s)
\end{lgathered}
\label{eq:reconstruction-input}
\end{equation}
where $\mathcal N_R(s)$ is a stencil
containing a fixed number of neighboring
subdomains, and $\chi_{s,j}^{(0)}$
denotes the initial reconstruction
feature associated with subdomain $K_j$.

At layer $q$, the convolutional features
are updated through
\begin{equation}
\begin{lgathered}
\chi_{s,j}^{(q+1)}=\operatorname{Mod}_R^{(q)}\left(\operatorname{ReLU}\left[\sum_{\delta\in\Lambda_R}W_{R,\delta}^{(q)}\widetilde\chi_{s,j+\delta}^{(q)}+b_R^{(q)}\right];\mathbf c\right)
\end{lgathered}
\label{eq:reconstruction-convolution}
\end{equation}
where $\Lambda_R$ is the set of
reconstruction-kernel offsets,
$W_{R,\delta}^{(q)}$ and
$b_R^{(q)}$ are the convolution
weights and bias, and
$\widetilde\chi$ denotes the
padded reconstruction features.

After $N_R$ convolution layers,
average pooling and linear projection
produce the Legendre coefficients:
\begin{equation}
\begin{lgathered}
\boldsymbol\xi_s^{n+1}=W_{R,\mathrm{out}}\left(\frac{1}{|\mathcal N_R(s)|}\sum_{j\in\mathcal N_R(s)}\chi_{s,j}^{(N_R)}\right)+b_{R,\mathrm{out}}
\end{lgathered}
\label{eq:legendre-coefficients}
\end{equation}
where $N_R$ is the number of
reconstruction convolution layers,
$W_{R,\mathrm{out}}$ and
$b_{R,\mathrm{out}}$ are the output
projection weights and bias, and
$\boldsymbol\xi_s^{n+1}
=(\xi_{s,1}^{n+1},\ldots,\xi_{s,P}^{n+1})$.
Here $P$ is the number of retained
Legendre modes, and
$\xi_{s,p}^{n+1}\in\mathbb R^{d+2}$
is the coefficient of the $p$th mode.

The local Legendre basis functions
are corrected to have zero mean:
\begin{equation}
\begin{lgathered}
\psi_{s,p}(x)=\operatorname{Leg}_p\left(\frac{x-a_s}{L}\right)-\frac{1}{V_s}\int_{K_s}\operatorname{Leg}_p\left(\frac{y-a_s}{L}\right)\,\mathrm dy
\end{lgathered}
\label{eq:legendre-basis}
\end{equation}
where $\operatorname{Leg}_p$ is the
$p$th multivariate Legendre mode,
$\psi_{s,p}$ is its mean-corrected
basis function on $K_s$, and $y$
is the integration variable.

The predicted flow field is reconstructed
within each local subdomain as
\begin{equation}
\begin{lgathered}
\widehat U^{n+1}(x)=\bar U_s^{n+1}+\sum_{p=1}^{P}\xi_{s,p}^{n+1}\psi_{s,p}(x),\qquad x\in K_s
\end{lgathered}
\label{eq:local-reconstruction}
\end{equation}
For discrete evaluation, the zero-mean
correction is applied using the same
quadrature as the local conservative states.
The reconstruction module yields
\begin{equation}
\begin{lgathered}
\widehat U^{n+1}=\mathcal R[\mathcal K_\Omega,\boldsymbol\mu](Z^n,\mathbf M^{n+1})
\end{lgathered}
\label{eq:reconstruction-module}
\end{equation}

Since
$\int_{K_s}\psi_{s,p}(x)\,\mathrm dx=0$,
the reconstructed field satisfies
\begin{equation}
\begin{lgathered}
\int_{K_s}\widehat U^{n+1}(x)\,\mathrm dx=M_s^{n+1}=M_s^n-\Delta t\sum_{\gamma\in\Gamma_s}\varepsilon_{s,\gamma}A_\gamma\widehat F_\gamma^n
\end{lgathered}
\label{eq:conservation-property}
\end{equation}
The reconstruction therefore preserves
the updated integral conservative states,
independently of the predicted Legendre
coefficients. The reconstructed flow is
re-encoded before the next time step.

\subsection{Boundary treatment}
\label{sec:boundary}

Boundary conditions are incorporated
into conservative evolution through
the interface fluxes. Let
$\Gamma_{\mathrm{int}}$,
$\Gamma_{\mathrm{per}}$,
$\Gamma_{\mathrm{open}}$, and
$\Gamma_{\mathrm{wall}}$
denote the sets of internal, periodically
paired, open-boundary, and wall interfaces,
respectively. The unified numerical
flux is defined by
\begin{equation}
\begin{lgathered}
\widehat F_\gamma^n=
\begin{cases}
\widehat F_{\gamma,\mathrm{learned}}^n, & \gamma\in\Gamma_{\mathrm{int}}\cup\Gamma_{\mathrm{per}},\\
\widehat F_{\gamma,\mathrm{open}}^n, & \gamma\in\Gamma_{\mathrm{open}},\\
\left(0,\left[(p_\gamma I-\tau_\gamma)n_\gamma\right]^{\mathsf T},0\right)^{\mathsf T}, & \gamma\in\Gamma_{\mathrm{wall}}.
\end{cases}
\end{lgathered}
\label{eq:unified-flux}
\end{equation}
Here $\widehat F_{\gamma,\mathrm{open}}^n$
is the open-boundary numerical flux,
$p_\gamma$ is the wall pressure,
$\tau_\gamma$ is the viscous stress tensor,
and $I$ is the identity matrix.
For physical boundary faces, $n_\gamma$
is taken as the outward unit normal,
so their orientation in
Eq.~\eqref{eq:fv-update} is positive.

For periodic boundaries, opposite
boundary interfaces are paired and
treated as internal interfaces.
Their convolutional stencils extend
periodically across the domain.
Each paired interface shares one
predicted flux with opposite orientations.

For open boundaries, the interior state
$U_\gamma^-$ is obtained from the adjacent
subdomain, while the exterior state
$U_\gamma^+$ is constructed from
the prescribed boundary data.
For characteristic boundary conditions,
incoming convective characteristic
information is prescribed, whereas
outgoing information is obtained
from the interior.

The open-boundary flux combines
the Rusanov convective flux with
the viscous contribution:
\begin{equation}
\begin{lgathered}
\widehat F_{\gamma,\mathrm{open}}^n=\frac{F^c(U_\gamma^-;\boldsymbol\mu)+F^c(U_\gamma^+;\boldsymbol\mu)}{2}n_\gamma-\frac{\lambda_\gamma^{\max}}{2}(U_\gamma^+-U_\gamma^-)-F^v\left(U_\gamma^b,(\nabla U)_\gamma^b;\boldsymbol\mu\right)n_\gamma
\end{lgathered}
\label{eq:open-boundary-flux}
\end{equation}
where $F^c$ and $F^v$ are the convective
and viscous fluxes defined in
Section~\ref{sec:problem},
$\lambda_\gamma^{\max}$ bounds
the normal convective wave speeds,
and $U_\gamma^b$ and
$(\nabla U)_\gamma^b$ denote
the boundary face state and its
discrete gradient, evaluated consistently
with the prescribed boundary conditions.

For a stationary impermeable, no-slip,
adiabatic wall, the velocity and thermal
boundary conditions are $\mathbf u=0$
and $\nabla T\cdot n_\gamma=0$,
where $\mathbf u$ and $T$ denote
velocity and temperature. The wall
flux in Eq.~\eqref{eq:unified-flux}
has zero mass and energy components.
Near nonperiodic boundaries, one-sided
or masked convolutional stencils are used.

\section{Experiments}
\label{sec:experiments}

The experiments assess four aspects of CLNO: prediction accuracy;
transfer across spatial resolutions, domain sizes, and boundary
conditions; generalization to unseen physical-parameter combinations;
and the contributions of conservative evolution and fine-scale
reconstruction. Sections~\ref{sec:one-d-results} and
\ref{sec:two-d-64} address the first two aspects using fixed-parameter
benchmarks. Section~\ref{sec:parameter-generalization} examines
parameter generalization, and Section~\ref{sec:ablation} evaluates
the contributions of the two components.

\subsection{Benchmarks and evaluation protocol}
\label{sec:experimental-setup}

For the fixed-parameter benchmarks, we evaluate the models on
one- and two-dimensional compressible Navier--Stokes flows.
The standard domains and grids are $[0,1]$ with 1024 cells and
$[0,1]^2$ with $64\times64$ cells, respectively. The training
data are generated with periodic boundary conditions, and the
nominal Reynolds numbers are $5\times10^3$ and $10^4$.
CLNO, LNO, FNO, and DeepONet are evaluated from the same initial
conservative state at the reference sampling times. All reported
forecasts are autoregressive: after initialization, each forecast
uses its own predictions as inputs to subsequent steps. Training data,
numerical solvers, model settings, and spatial-transfer protocols are specified in
Appendix~\ref{app:experimental-details}.

For $T$ forecast frames, we report the relative trajectory error
\begin{equation}
e_{L^2}=\frac{\|\widehat U^{1:T}-U^{1:T}\|_2}{\|U^{1:T}\|_2},
\label{eq:trajectory-error}
\end{equation}
where the norm covers time, space, and all conservative variables.
Let $J^n=\sum_{x\in\Omega_f^n}U^n(x)\,\Delta V$ denote the
spatial integral of the conservative state over the fluid region
at frame $n$.
The conservation error compares the predicted and reference net
changes between the initial and final frames,
\begin{equation}
e_{\mathrm{cons}}=\frac{\|(\widehat J^T-\widehat J^0)-(J^T-J^0)\|_2}{\|J^0\|_2}.
\label{eq:conservation-error}
\end{equation}
The norm includes mass, each momentum component, and total energy.
This definition accounts for the reference net change when the
fluid region has nonperiodic boundaries or an obstacle. Tables
express both errors as percentages.
Inference time covers the complete rollout, including boundary
processing where applicable. In the one-dimensional fixed-parameter
and Reynolds--Mach comparisons, finite forecasts are scored even
when density or pressure becomes nonpositive. A forecast containing
NaN or infinity before the final frame is reported as \texttt{NaN};
no full-trajectory error is assigned.

\FloatBarrier
\subsection{One-dimensional compressible Navier--Stokes equations}
\label{sec:one-d-results}

The one-dimensional evaluation spans periodic multiscale waves,
transient reflecting obstacles, and advected contact discontinuities
on the training grid, a refined grid, and an extended domain. A
pressure-discontinuity Riemann problem probes a distinct nonperiodic
initial state. Network weights are fixed across the three spatial
configurations. Table~\ref{tab:one-d-results} summarizes complete
100-frame forecasts; the reference discretization and data splits
are specified in Appendix~\ref{app:one-d-protocol}.

The periodic-wave and obstacle experiments use the same
low-frequency Fourier initial field. The outer boundaries are
periodic in both tests; the obstacle introduces an internal
reflecting boundary during its active interval.
Appendix~\ref{app:one-d-protocol} gives the initial-field
definition and random seed.

\begin{table}[!htbp]
\centering
\caption{One-dimensional Navier--Stokes equation forecasts over
100 frames. $(N,L)$ denotes grid cells and domain length.
Relative $L^2$ errors are percentages and inference times are
seconds. \texttt{NaN} indicates that a state became nonfinite before
the final frame; partial runtimes are omitted. Finite forecasts are
scored regardless of state positivity. Boldface marks the lowest
finite value in each row and metric.}
\label{tab:one-d-results}
\small
\setlength{\tabcolsep}{3pt}
\renewcommand{\arraystretch}{1.12}
\resizebox{\linewidth}{!}{
\begin{tabular}{@{}ll*{4}{c}*{4}{c}*{4}{c}@{}}
\toprule
Flow configuration
& $(N,L)$
& \multicolumn{4}{c}{Relative $L^2$ error (\%)}
& \multicolumn{4}{c}{Inference time (s)}
& \multicolumn{4}{c}{Conservation error (\%)}\\
\cmidrule(lr){3-6}
\cmidrule(lr){7-10}\cmidrule(l){11-14}
& &
FNO & DeepONet & LNO & CLNO &
FNO & DeepONet & LNO & CLNO & FNO & DeepONet & LNO & CLNO\\
\midrule
Periodic multiscale wave
& $(1024,1)$
& 6.392 & 4.594 & 7.563 & \textbf{1.840}
& 0.1893 & \textbf{0.1190} & 0.1499 & 0.6432 & 1.308 & 0.054 & 4.690 & \textbf{1.59e-6}\\
& $(1536,1)$
& 6.663 & 4.839 & 27.989 & \textbf{1.346}
& 0.2101 & \textbf{0.1114} & 0.1302 & 0.6436 & 1.310 & 0.053 & 4.974 & \textbf{2.33e-6}\\
& $(1536,1.5)$
& 29.267 & 28.885 & 6.939 & \textbf{1.220}
& 0.2047 & \textbf{0.1320} & 0.1344 & 0.6342 & 1.310 & 0.053 & 4.974 & \textbf{4.55e-6}\\
\addlinespace
Transient reflecting obstacle
& $(1024,1)$
& \texttt{NaN} & 16.960 & 11.050 & \textbf{3.045}
& --- & 0.3750 & \textbf{0.1759} & 2.2425 & \texttt{NaN} & 4.637 & 5.269 & \textbf{0.365}\\
& $(1536,1)$
& \texttt{NaN} & 17.154 & 26.103 & \textbf{2.993}
& --- & 0.2186 & \textbf{0.1548} & 2.2206 & \texttt{NaN} & 4.606 & 6.154 & \textbf{0.353}\\
& $(1536,1.5)$
& \texttt{NaN} & 22.880 & 9.446 & \textbf{2.340}
& --- & 0.2304 & \textbf{0.1552} & 2.2240 & \texttt{NaN} & 5.291 & 5.431 & \textbf{0.097}\\
\addlinespace
Advected contact discontinuity
& $(1024,1)$
& 21.567 & 5916.028 & 21.124 & \textbf{0.880}
& 0.1883 & \textbf{0.0857} & 0.1303 & 1.3788 & 19.978 & 3111.598 & 22.926 & \textbf{0.046}\\
& $(1536,1)$
& 24.429 & 8322.915 & 19.713 & \textbf{0.741}
& 0.1931 & \textbf{0.0840} & 0.1347 & 1.3770 & 20.407 & 4642.523 & 18.007 & \textbf{0.046}\\
& $(1536,1.5)$
& 24.014 & 8353.556 & 19.295 & \textbf{0.703}
& 0.1839 & \textbf{0.0737} & 0.1300 & 1.3713 & 19.496 & 4643.043 & 18.065 & \textbf{0.029}\\
\addlinespace
Pressure-jump Riemann problem
& $(1024,1)$
& 42.554 & 45.968 & 33.732 & \textbf{9.383}
& 0.1658 & \textbf{0.0632} & 0.1845 & 1.5396 & 3.997 & 24.015 & 40.965 & \textbf{1.547}\\
\bottomrule
\end{tabular}}
\end{table}

\paragraph{Periodic and spatial-transfer tests.}

A held-out low-frequency Fourier realization is propagated with
periodic boundaries. Grid refinement changes the resolution at
fixed domain length, whereas domain extension retains the grid
spacing. CLNO has the lowest trajectory error among the four
models in all three configurations. On the refined grid, its
error is $1.346\%$, compared with $27.989\%$ for LNO. These
results support the use of the same trained CLNO model across
the three tested spatial configurations. Additional space--time
and profile comparisons are shown in
Figures~\ref{fig:one-d-fourier-fields} and
\ref{fig:one-d-fourier-snapshots} in
Appendix~\ref{app:one-d-protocol}.

\paragraph{Nonperiodic discontinuities.}

\begin{figure}[H]
\centering
\includegraphics[width=0.72\linewidth]{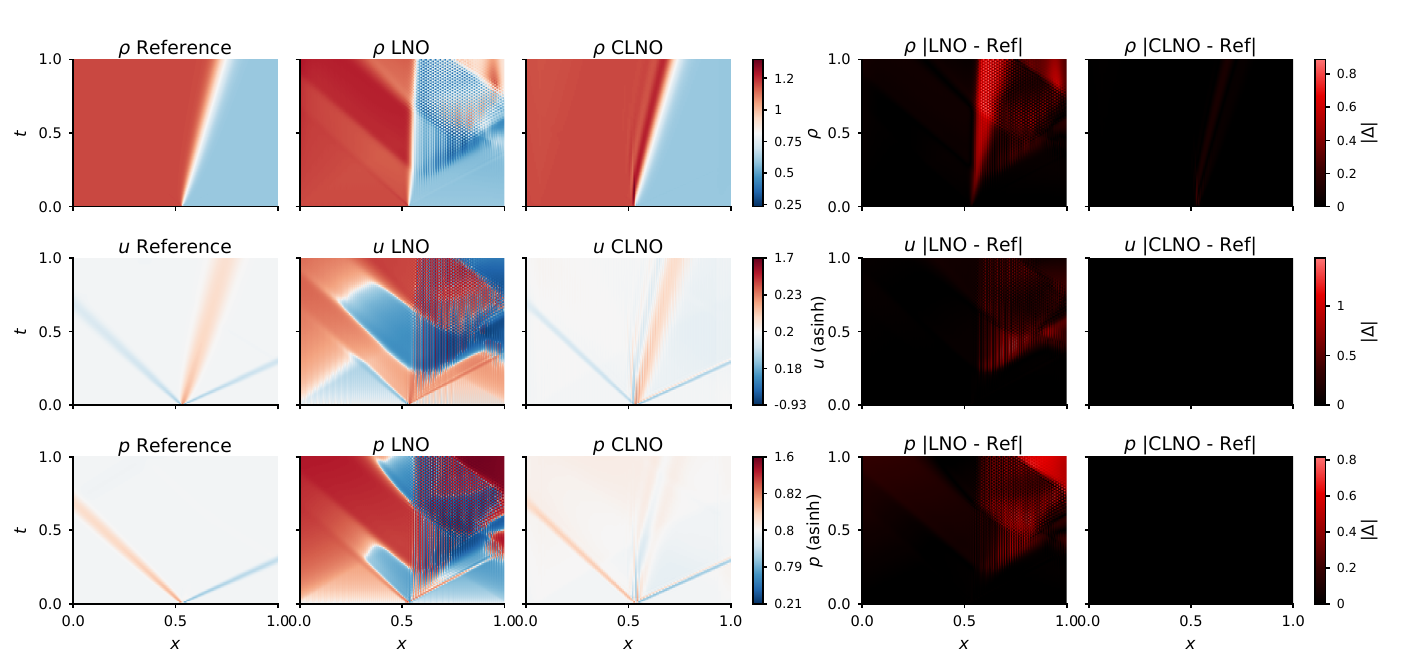}
\caption{Space--time density, velocity, and pressure for
nonperiodic contact transport. Columns show the reference,
LNO, CLNO, and their absolute errors. The velocity and pressure
field colorbars use a centered inverse-hyperbolic-sine scale to
resolve variations near their nearly constant reference values;
tick labels give the original physical values.}
\label{fig:one-d-contact-fields}
\end{figure}

\begin{figure}[H]
\centering
\includegraphics[width=0.72\linewidth]{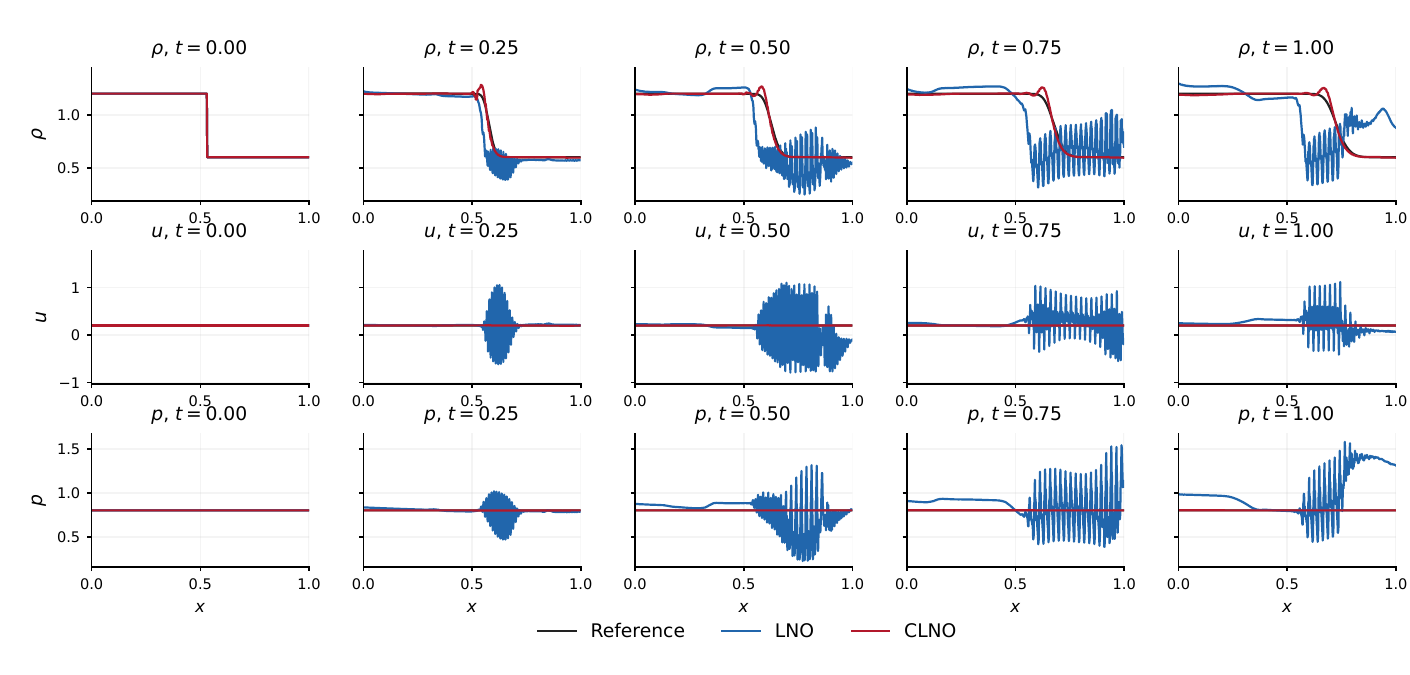}
\caption{Contact-discontinuity profiles at $t=0$, $0.25$,
$0.5$, $0.75$, and $1$.}
\label{fig:one-d-contact-snapshots}
\end{figure}

The contact-discontinuity test starts from a density discontinuity
with equal pressure and velocity on both sides. Constant-extrapolation
outer boundaries isolate contact transport. The primitive initial
state and exterior ghost-cell condition are
\begin{equation}
(\rho_0,u_0,p_0)=
(0.6+0.6\mathbf 1_{\{x/L<0.53\}},\,0.2,\,0.8),\qquad (U_{g,L},U_{g,R})=(U_1,U_N).
\label{eq:contact-initial-main}
\end{equation}
Here \(\mathbf 1_{\{x/L<0.53\}}\) is the indicator of the left state.
The contact problem probes the effect of nonperiodic transport
over long autoregressive forecasts. CLNO maintains a trajectory
error below $1\%$ in all three spatial configurations, whereas
LNO exceeds $19\%$. On the training grid, FNO and DeepONet have
errors of $21.567\%$ and $5916.028\%$, respectively.
Figure~\ref{fig:one-d-contact-fields}
shows that CLNO preserves the localization of the density jump
and the nearly uniform velocity and pressure. LNO develops
oscillations. Profiles at five times are shown in
Figure~\ref{fig:one-d-contact-snapshots}.

\FloatBarrier

\FloatBarrier

\paragraph{Transient reflecting boundary.}

A temporary reflecting obstacle interrupts the periodic flow
and is subsequently removed, testing both wave reflection and
propagation after the fluid region reconnects. The initial state
is the periodic-wave state defined in
Eq.~\eqref{app:eq:one-d-periodic-state}. During
\(0.2\leq t<0.8\), the solid occupies
\(x/L\in[0.4375,0.5625]\); its wall ghost state is
\(U_g=(\rho,-\rho u,E)_f\).
The exterior boundary remains periodic. Outside this time interval
the domain is entirely fluid.

\begin{figure}[H]
\centering
\includegraphics[width=0.72\linewidth]{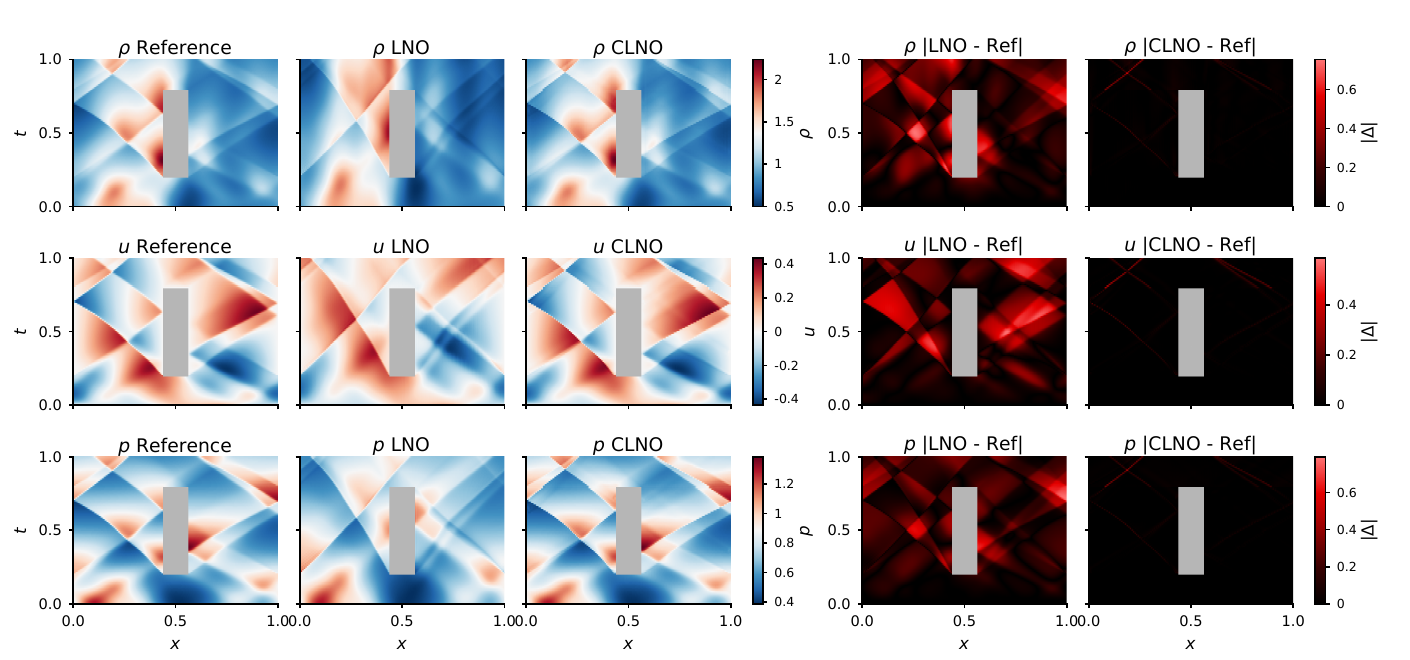}
\caption{Space--time prediction with a temporary reflecting
obstacle on the refined 1536-cell grid. The solid interval is
shaded gray during its active period.}
\label{fig:one-d-obstacle-cross-fields}
\end{figure}

\begin{figure}[H]
\centering
\includegraphics[width=0.72\linewidth]{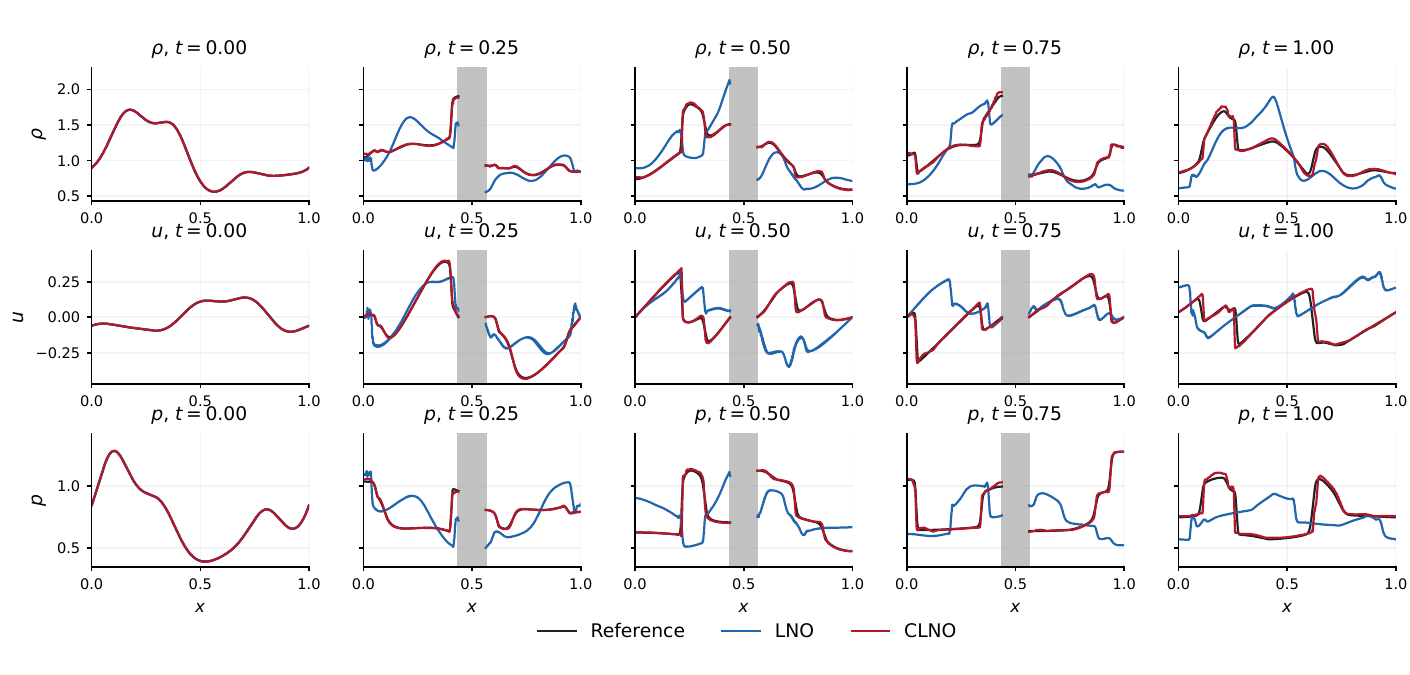}
\caption{Density, velocity, and pressure profiles of the
refined-grid obstacle test at $t=0$, $0.25$, $0.5$, $0.75$,
and $1$.}
\label{fig:one-d-obstacle-cross-snapshots}
\end{figure}

CLNO yields lower trajectory errors than LNO and DeepONet
in all three configurations. FNO produces nonfinite values
before the final frame in each configuration. On the refined grid,
CLNO and LNO errors are $2.993\%$ and $26.103\%$, respectively.
Figures~\ref{fig:one-d-obstacle-cross-fields} and
\ref{fig:one-d-obstacle-cross-snapshots} show the resulting
space--time fields and profiles. The longer CLNO inference
time includes local reconstruction and obstacle-boundary
processing. The standard-grid comparison is shown in
Figure~\ref{fig:one-d-obstacle-fields} in
Appendix~\ref{app:one-d-protocol}.

\FloatBarrier

\FloatBarrier

\paragraph{Independent pressure-jump test.}

A complementary test introduces an initial pressure discontinuity
from $1$ to $0.2$ at $x=0.5$. It uses the same
constant-extrapolation boundary condition as
Eq.~\eqref{eq:contact-initial-main}, with initial state
\begin{equation}
(\rho_0,u_0,p_0)=(1,\,0,\,0.2+0.8\mathbf 1_{\{x<0.5\}}).
\label{eq:pressure-jump-initial-main}
\end{equation}
CLNO also outperforms LNO
on this standard-grid forecast (Table~\ref{tab:one-d-results}).
Field and profile comparisons for this test are shown in
Figures~\ref{fig:one-d-pressure-fields} and
\ref{fig:one-d-pressure-snapshots} in
Appendix~\ref{app:one-d-protocol}. Across the original 28
one-dimensional tests, CLNO has lower full-trajectory error
than LNO.

\FloatBarrier
\subsection{Two-dimensional compressible Navier--Stokes equations}
\label{sec:two-d-64}

The two-dimensional Navier--Stokes evaluation comprises held-out
vortical trajectories, an independent high-wavenumber perturbation,
and a matched spatial-transfer experiment. For the last experiment,
one initial field is evaluated on the training grid, a refined
grid, and an extended domain with unchanged weights and equation
parameters. Table~\ref{tab:two-d-64-results} specifies the grids;
the latent-grid sizes are given in Appendix~\ref{app:two-d-protocol}.

The held-out and spatial-transfer cases start from periodic
vortical fields with smooth density and pressure perturbations.
The independent high-wavenumber test increases the maximum
Fourier wavenumber from four in the training perturbations to
eight, testing prediction with finer spatial variations.
All cases use periodic boundary conditions in both directions.
Appendix~\ref{app:two-d-protocol} specifies the initial-field
definitions, parameter ranges, and random seeds.

\begin{figure}[H]
\centering
\includegraphics[width=0.72\linewidth]{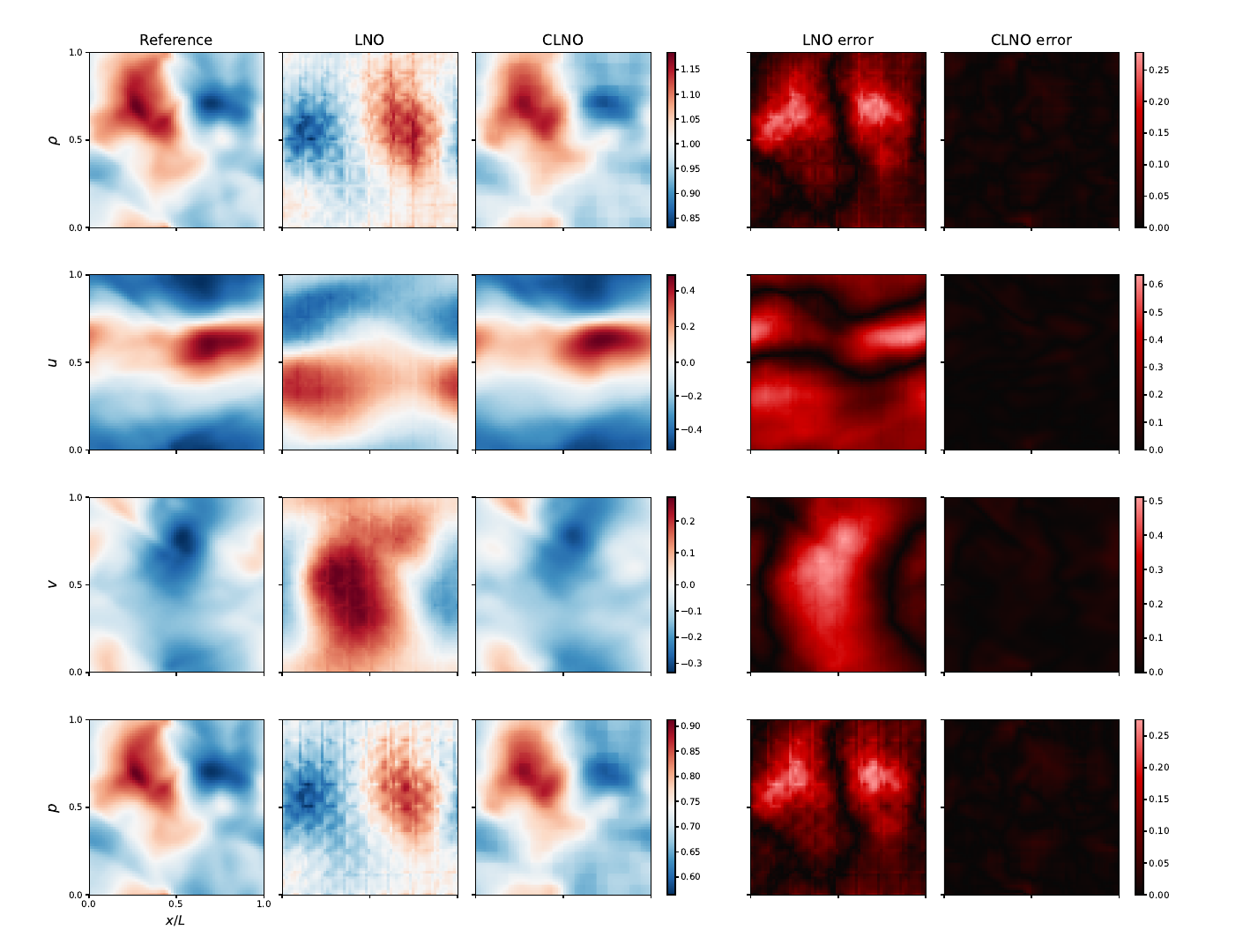}
\caption{Held-out two-dimensional vortical flow at the final
forecast time $t=8$ on the $64\times64$ grid. Columns show
the reference, LNO, CLNO, and their absolute errors; rows show
density, both velocity components, and pressure.}
\label{fig:two-d-64-heldout}
\end{figure}

Table~\ref{tab:two-d-64-results} reports full-trajectory errors
and rollout times. CLNO has the lowest mean error on the held-out
vortical flows and the lowest error on the high-wavenumber
perturbation. On the former, its mean error is $3.306\%$,
compared with $15.901\%$ for LNO. In the representative flow
shown in Fig.~\ref{fig:two-d-64-heldout}, their respective
instantaneous conservative-state errors at $t=8$ are $2.382\%$
and $22.163\%$.

\begin{table}[t]
\centering
\caption{Representative two-dimensional Navier--Stokes equation
forecasts. Relative $L^2$ errors are computed over each complete
space--time trajectory. The held-out row reports the mean over ten
initial fields; the final three rows share one initial field.
$(N,L)$ gives the grid size per
direction and domain side length; $T$ is the final physical time.
Boldface marks the lowest finite value in each row and metric.}
\label{tab:two-d-64-results}
\small
\setlength{\tabcolsep}{2pt}
\renewcommand{\arraystretch}{1.12}
\resizebox{\linewidth}{!}{
\begin{tabular}{@{}l*{12}{c}@{}}
\toprule
Flow configuration
& \multicolumn{4}{c}{Relative $L^2$ error (\%)}
& \multicolumn{4}{c}{Inference time (s)}
& \multicolumn{4}{c}{Conservation error (\%)}\\
\cmidrule(lr){2-5}
\cmidrule(lr){6-9}\cmidrule(l){10-13}
& FNO & DeepONet & LNO & CLNO
& FNO & DeepONet & LNO & CLNO & FNO & DeepONet & LNO & CLNO\\
\midrule
\shortstack[l]{Held-out vortical flows\\$(64,1)$, $T=8$}
& 4.946 & 14.124 & 15.901 & \textbf{3.306}
& 0.3066 & \textbf{0.1747} & 0.2858 & 0.9649 & 0.476 & 0.608 & 4.172 & \textbf{2.39e-6}\\
\shortstack[l]{High-wavenumber perturbation\\$(64,1)$, $T=1$}
& 8.529 & 18.855 & 12.005 & \textbf{7.101}
& 0.0363 & \textbf{0.0244} & 0.0352 & 0.1211 & 0.754 & 0.399 & 1.773 & \textbf{1.20e-6}\\
\addlinespace
\shortstack[l]{Vortical flow, standard grid\\$(64,1)$, $T=1$}
& 1.428 & 10.354 & 5.789 & \textbf{1.193}
& 0.0419 & \textbf{0.0154} & 0.0462 & 0.1580 & 0.293 & 0.207 & 1.704 & \textbf{9.29e-7}\\
\shortstack[l]{Vortical flow, refined grid\\$(96,1)$, $T=1$}
& \textbf{1.551} & 10.660 & 13.895 & 1.560
& 0.0420 & \textbf{0.0153} & 0.0452 & 0.1362 & 0.287 & 0.214 & 2.717 & \textbf{5.17e-7}\\
\shortstack[l]{Vortical flow, expanded domain\\$(96,1.5)$, $T=1$}
& 14.842 & 14.439 & 9.035 & \textbf{3.328}
& 0.0365 & \textbf{0.0135} & 0.0393 & 0.1364 & 0.287 & 0.214 & 2.718 & \textbf{2.59e-7}\\
\bottomrule
\end{tabular}}
\end{table}

For the matched realization, CLNO yields lower errors than LNO
in all three configurations, with the largest difference on the
refined grid. FNO and CLNO yield comparable errors on the
refined grid ($1.551\%$ and $1.560\%$, respectively), while CLNO
has the lowest error among the four methods in the other two
matched settings. CLNO's error increases from $1.193\%$ on the
standard grid to $1.560\%$ on the refined grid and $3.328\%$ on
the extended domain. These results support spatial transfer for
this matched initial field without retraining, while indicating
reduced accuracy on the enlarged domain.

\begin{figure}[H]
\centering
\includegraphics[width=0.72\linewidth]{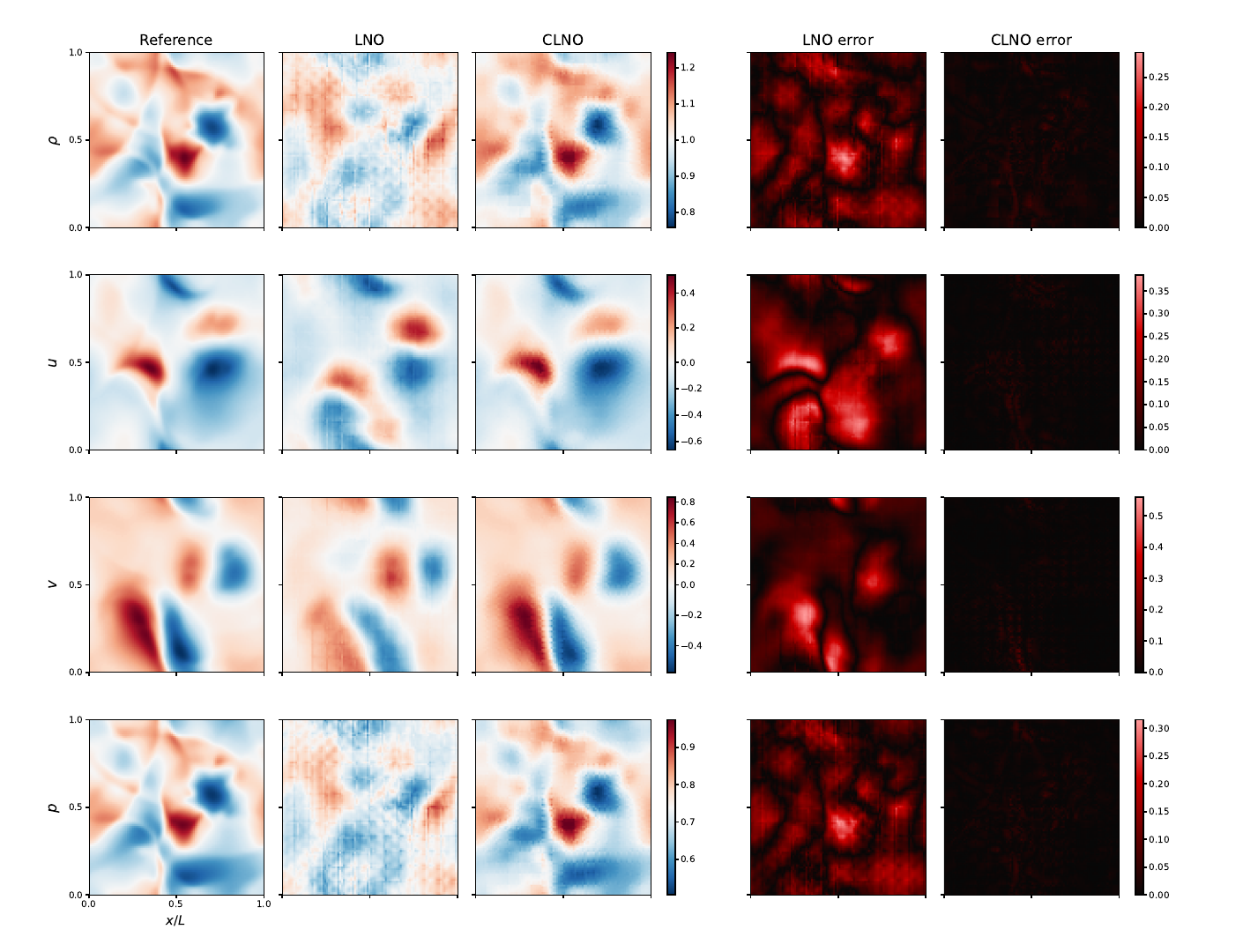}
\caption{The matched vortical realization on a refined
$96\times96$ grid over $[0,1]^2$ at $t=0.75$.}
\label{fig:two-d-64-refined}
\end{figure}

\begin{figure}[H]
\centering
\includegraphics[width=0.72\linewidth]{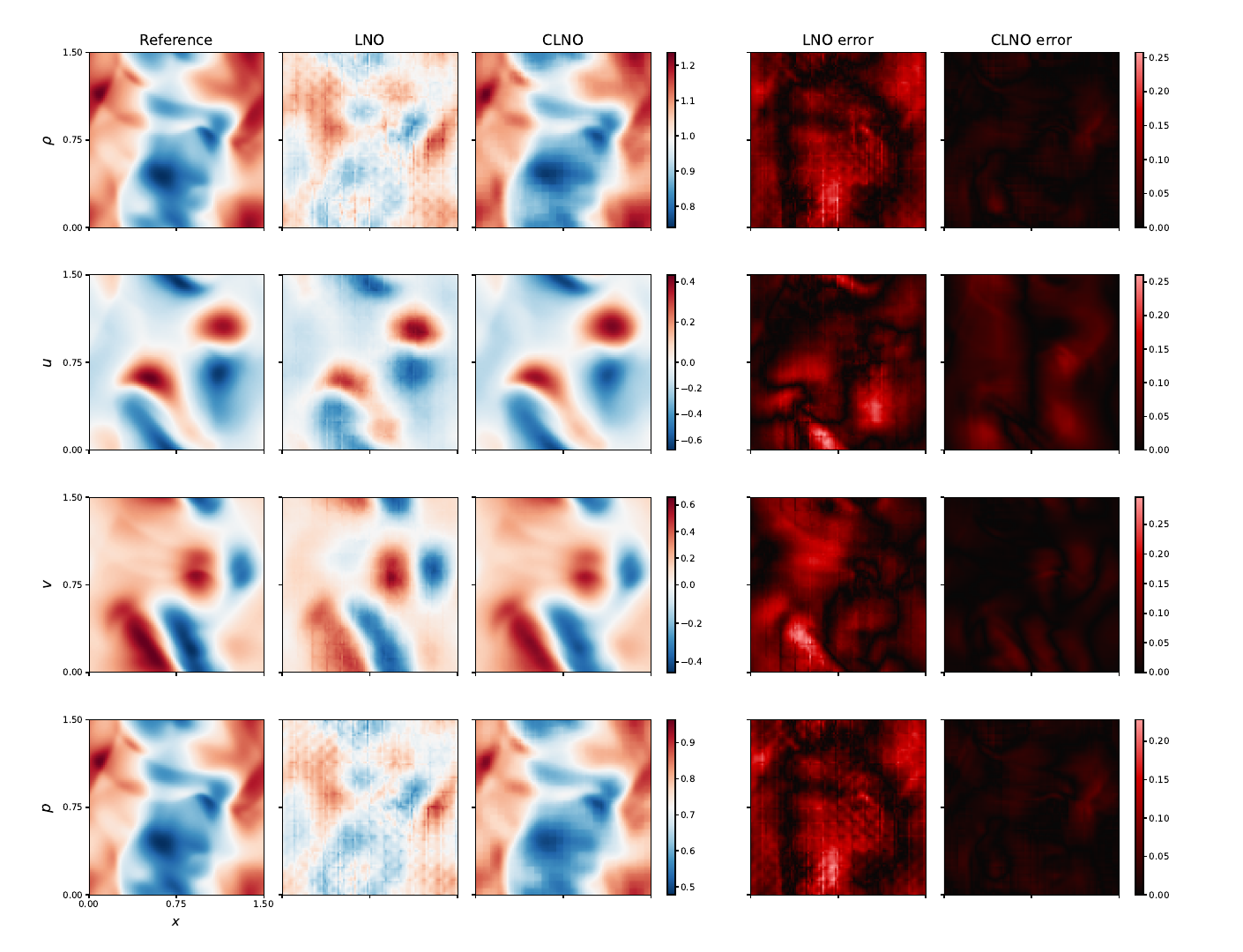}
\caption{The matched vortical realization over the expanded
domain $[0,1.5]^2$ on a $96\times96$ grid at $t=0.75$.}
\label{fig:two-d-64-expanded}
\end{figure}

Figures~\ref{fig:two-d-64-refined} and
\ref{fig:two-d-64-expanded} show the matched realization at
$t=0.75$ after grid refinement and domain extension,
respectively. At this frame, the LNO-to-CLNO error ratios are
$8.74$ and $3.45$. The standard-grid comparison at the same
time is shown in Figure~\ref{fig:two-d-64-standard} in
Appendix~\ref{app:two-d-protocol}.
All three figures compare density, both velocity components,
pressure, and absolute errors using the same panel arrangement.

\FloatBarrier
\FloatBarrier

\FloatBarrier

\FloatBarrier

\FloatBarrier
\subsection{Generalization across Reynolds and Mach numbers}
\label{sec:parameter-generalization}

This experiment tests whether a trained CLNO model can predict
flows at physical-parameter combinations excluded from training.
The benchmark uses a fixed periodic one-dimensional grid with
prescribed dimensionless Reynolds and Mach numbers. Of the ten
parameter pairs in the experimental design,
$(\mathrm{Re},\mathrm{Ma})=(7500,0.45)$ is held out, and the
remaining nine pairs provide the training and validation data.
We also evaluate $(5000,0.50)$, which lies outside the training
parameter range.

The benchmark uses its own reference dataset and separately
trained model weights. The reference trajectories are generated
with a finite-volume Rusanov solver, independently of the
fixed-parameter benchmarks. The parameter pairs use the same ten
primitive initial fields. The initial-field definitions, parameter
grid, and data splits are specified in
Appendix~\ref{app:parameter-protocol}.

\begin{table}[!htbp]
\centering
\caption{One-dimensional Reynolds--Mach parameter sweep on a periodic 1024-cell grid. The first eleven rows use the same paired initial field: a validation initial field at seen parameter pairs, a test initial field at the held-out pair $(7500,0.45)$, and the corresponding independent initial field at the outside-range pair $(5000,0.50)$. The final two rows report means over ten test initial fields at the two unseen pairs. All models predict 64 future frames. Errors are in percent; time is per complete trajectory in seconds. Boldface marks the lowest finite value in each row and metric.}
\label{tab:re-ma-results}
\small
\setlength{\tabcolsep}{2pt}
\renewcommand{\arraystretch}{1.10}
\resizebox{\linewidth}{!}{
\begin{tabular}{@{}rc*{12}{c}@{}}
\toprule
$\mathrm{Re}$ & $\mathrm{Ma}$ & \multicolumn{4}{c}{Relative $L^2$ error (\%)} & \multicolumn{4}{c}{Inference time (s)} & \multicolumn{4}{c}{Conservation error (\%)}\\
\cmidrule(lr){3-6}\cmidrule(lr){7-10}\cmidrule(l){11-14}
& & FNO & DeepONet & LNO & CLNO & FNO & DeepONet & LNO & CLNO & FNO & DeepONet & LNO & CLNO\\
\midrule
5500 & 0.35 & 0.309 & 1.089 & 94.809 & \textbf{0.066} & 0.1219 & \textbf{0.0556} & 0.0882 & 0.7106 & 0.169 & 0.055 & 224.594 & \textbf{8.24e-7}\\
5500 & 0.45 & 0.237 & 0.909 & 12.994 & \textbf{0.077} & 0.1222 & \textbf{0.0559} & 0.0882 & 0.7108 & 0.115 & 0.033 & 30.610 & \textbf{8.24e-7}\\
6500 & 0.35 & 0.320 & 1.097 & 159.986 & \textbf{0.064} & 0.1220 & \textbf{0.0558} & 0.0882 & 0.7104 & 0.134 & 0.065 & 476.535 & \textbf{6.66e-7}\\
6500 & 0.45 & 0.234 & 0.912 & 41.951 & \textbf{0.078} & 0.1225 & \textbf{0.0556} & 0.0879 & 0.7096 & 0.061 & 0.032 & 97.441 & \textbf{4.24e-6}\\
7500 & 0.35 & 0.330 & 1.100 & 129.474 & \textbf{0.063} & 0.1226 & \textbf{0.0556} & 0.0883 & 0.7106 & 0.116 & 0.067 & 335.866 & \textbf{1.80e-6}\\
7500 & 0.45 & 0.239 & 0.910 & 58.724 & \textbf{0.078} & 0.1225 & \textbf{0.0556} & 0.0884 & 0.7102 & 0.049 & 0.033 & 137.380 & \textbf{4.30e-6}\\
8500 & 0.35 & 0.323 & 1.100 & 120.563 & \textbf{0.063} & 0.1224 & \textbf{0.0556} & 0.0882 & 0.7110 & 0.123 & 0.064 & 281.223 & \textbf{3.53e-6}\\
8500 & 0.45 & 0.241 & 0.904 & 49.043 & \textbf{0.077} & 0.1058 & \textbf{0.0556} & 0.0881 & 0.7112 & 0.073 & 0.034 & 123.411 & \textbf{1.90e-6}\\
9500 & 0.35 & 0.311 & 1.098 & 117.006 & \textbf{0.065} & 0.1058 & \textbf{0.0556} & 0.0883 & 0.8164 & 0.156 & 0.058 & 291.936 & \textbf{2.41e-6}\\
9500 & 0.45 & 0.249 & 0.898 & 117.913 & \textbf{0.075} & 0.1058 & \textbf{0.0556} & 0.0883 & 0.7001 & 0.129 & 0.035 & 318.796 & \textbf{3.10e-6}\\
5000 & 0.50 & 2.597 & 4.819 & 12223.604 & \textbf{0.299} & 0.1057 & \textbf{0.0556} & 0.0882 & 0.6997 & 0.995 & 0.057 & 36649.026 & \textbf{1.29e-6}\\
\midrule
\multicolumn{14}{@{}l}{Mean over ten test initial fields at unseen parameter pairs}\\
7500 & 0.45 & 0.747 & 1.225 & 436.602 & \textbf{0.140} & 0.1207 & \textbf{0.0633} & 0.1035 & 0.7414 & 0.158 & 0.036 & 1080.604 & \textbf{2.79e-6}\\
5000 & 0.50 & 3.970 & 6.790 & 38009.355 & \textbf{0.511} & 0.1046 & \textbf{0.0533} & 0.0865 & 0.7436 & 0.917 & 0.058 & 113680.415 & \textbf{2.06e-6}\\
\bottomrule
\end{tabular}}
\end{table}

Across ten initial fields, CLNO achieves mean relative trajectory
errors of $0.140\%$ at the held-out pair $(7500,0.45)$ and
$0.511\%$ at the outside-range pair $(5000,0.50)$ over 64 forecast
frames (Table~\ref{tab:re-ma-results}). These errors are lower than
those of all three comparison models at both parameter pairs.
The outside-range case has a higher error, but the mean remains
below $1\%$ in both settings. These results support prediction
at the two evaluated unseen parameter combinations without
retraining.

The conditioned LNO produces finite values throughout the
64-frame forecasts, but its error grows sharply as the forecast
horizon extends beyond the shorter horizons used during training.
Its full-trajectory errors are therefore included in
Table~\ref{tab:re-ma-results}.

\FloatBarrier
\subsection{Ablation studies}
\label{sec:ablation}

We compare conservative flux evolution and fine-scale
reconstruction on the fixed-parameter $64\times64$
two-dimensional Navier--Stokes benchmark.
The ablated models and their matched training protocol are specified
in Appendix~\ref{app:parameter-protocol}.
\begin{table}[!htbp]
\centering
\caption{Two-dimensional $64\times64$ Navier--Stokes equation ablation on two representative 160-frame held-out trajectories. The trajectories maximize the smaller of the two ablation-to-full error ratios among ten complete forecasts. Values are means over these two trajectories. Boldface marks the lowest value in each metric.}
\label{tab:ablation-fixed2d}
\small
\setlength{\tabcolsep}{3pt}
\begin{tabular}{@{}lccc@{}}
\toprule
Model & Relative $L^2$ error (\%) & Inference time (s) & Conservation error (\%)\\
\midrule
CLNO (w/ both components) & \textbf{1.384} & 1.050 & \textbf{2.13e-6}\\
CLNO (w/o conservative flux) & 1.953 & 0.962 & 0.337\\
CLNO (w/o fine reconstruction) & 3.351 & \textbf{0.657} & 3.14e-6\\
\bottomrule
\end{tabular}
\end{table}
On the two selected held-out trajectories, the full model has the
lowest error. Replacing shared interface-flux evolution with a
per-cell residual increases trajectory error from $1.384\%$ to
$1.953\%$ and raises conservation error to $0.337\%$; suppressing
fine-scale reconstruction increases trajectory error to $3.351\%$
while retaining a conservative latent update. The latter is faster
because it omits the fine-grid decoder.

\FloatBarrier

\section{Discussion and conclusion}
\label{sec:conclusion}

We proposed a conservative local neural operator (CLNO)
for accurate and generalizable prediction of compressible
flows. By combining shared interface-flux prediction,
finite-volume evolution, and zero-mean Legendre
reconstruction, CLNO preserves the updated local
conservative states while recovering fine-scale flow
structures. Numerical experiments demonstrate improved
prediction accuracy and generalization compared with
LNO across different spatial resolutions, computational
domains, and boundary conditions, together with improved
long-time prediction performance. Furthermore,
Reynolds--Mach number conditioning enables a unified
model to represent parameter-dependent flow evolution
and directly predict unseen parameter combinations
without retraining.

Several directions can further improve the adaptability
and physical consistency of CLNO. First, the physical
size of local subdomains could be made differentiable
and adaptively determined by the Reynolds and Mach
numbers, allowing the local representation to adjust
to different flow regimes. Second, a parameter-dependent
gating mechanism could dynamically control the number
of active Legendre modes, balancing fine-scale
reconstruction accuracy and computational efficiency.
Third, although CLNO preserves the integral conservation
of each local subdomain, the reconstructed field is not
explicitly constrained to satisfy the governing equations
pointwise. Incorporating local physics-informed losses
or differentiable physical corrections during training
or inference could further improve physical consistency
while retaining the established conservation property.

Overall, CLNO provides a conservative and generalizable
framework for compressible flow prediction. Its combination
of local conservative evolution, fine-scale reconstruction,
and physical-parameter modulation establishes a foundation
for further developments in adaptive and physically
consistent neural operators.

\FloatBarrier
\bibliographystyle{cas-model2-names}
\bibliography{cas-refs}

\appendix
\section{Experimental data and implementation details}
\label{app:experimental-details}

This appendix records the numerical reference data, model configurations,
and rollout protocols underlying Section~\ref{sec:experiments}. The
fixed-parameter and Reynolds--Mach experiments use separate datasets and
checkpoints. Reynolds and Mach numbers in the parameter study are
prescribed dimensionless equation parameters; they are not recomputed
from an initial field or a changed computational domain.

\subsection{Fixed-parameter one-dimensional data}
\label{app:one-d-protocol}

The compressible Navier--Stokes equations are solved on the periodic
training domain $[0,1]$ with $N=1024$ cells, $\gamma=1.4$, and viscosity
$2\times10^{-4}$ (nominal $\mathrm{Re}=5000$). The training reference solver uses
fifth-order WENO reconstruction, Lax--Friedrichs flux splitting,
physical viscous stress and Fourier heat conduction, followed by
third-order SSP Runge--Kutta integration. The archive contains 124
trajectories of 201 frames over $t\in[0,2]$, sampled every $0.01$.
Trajectory-level splits contain 108 training, eight validation, and
eight held-out sequences. The primitive and conservative archives each
have shape $(124,201,3,1024)$.

The independent obstacle, contact, and pressure-jump test references
use a finite-volume Rusanov flux, centered Navier--Stokes diffusion,
and third-order SSP Runge--Kutta time stepping at the same viscosity.
They are sampled at $\Delta t=0.01$, with CFL and viscous substeps
inside each stored-frame interval.

Training initial states include low-frequency random Fourier fields
with modes $1$--$4$, smooth nonlinear waves, and 24 additional
high-frequency Fourier fields with wave numbers $6$--$16$. The three
test geometries in Table~\ref{tab:one-d-results} use
$(N,L)=(1024,1),(1536,1),(1536,1.5)$. Grid refinement changes $N$
at fixed $L$, whereas domain extension retains the original grid
spacing. The CLNO latent-cell width is $1/128$, giving 128, 128, and
192 cells, respectively, without changing learned weights.

The periodic-wave and obstacle experiments use the same Fourier
initial field. With \(\widetilde x=x/L\), define the unit-standard-deviation
profiles
\begin{equation}
\mathcal W_j(\widetilde x)=\frac{1}{\sigma_j}\sum_{k=1}^{4}k^{-5/4}\chi_j(2\pi k\widetilde x+\phi_{jk}),
\label{eq:one-d-fourier-initial}
\end{equation}
where \(\chi_\rho=\chi_p=\cos\), \(\chi_u=\sin\), and the
independent phases \(\phi_{jk}\) specify the held-out realization;
\(\sigma_j\) is the standard deviation of the unnormalized sum.
The primitive initial state is
\begin{equation}
(\rho_0,u_0,p_0)
=(e^{0.35\mathcal W_\rho},\,0.0875\mathcal W_u,\,e^{0.35\mathcal W_p}/1.4).
\label{app:eq:one-d-periodic-state}
\end{equation}
The outer boundary is periodic:
\(U(0,t)=U(L,t)\) and
\(\partial_xU(0,t)=\partial_xU(L,t)\).
The plotted realization uses the phases generated with seed
20261026 in Eq.~\eqref{eq:one-d-fourier-initial}.

\begin{figure}[H]
\centering
\includegraphics[width=0.84\linewidth]{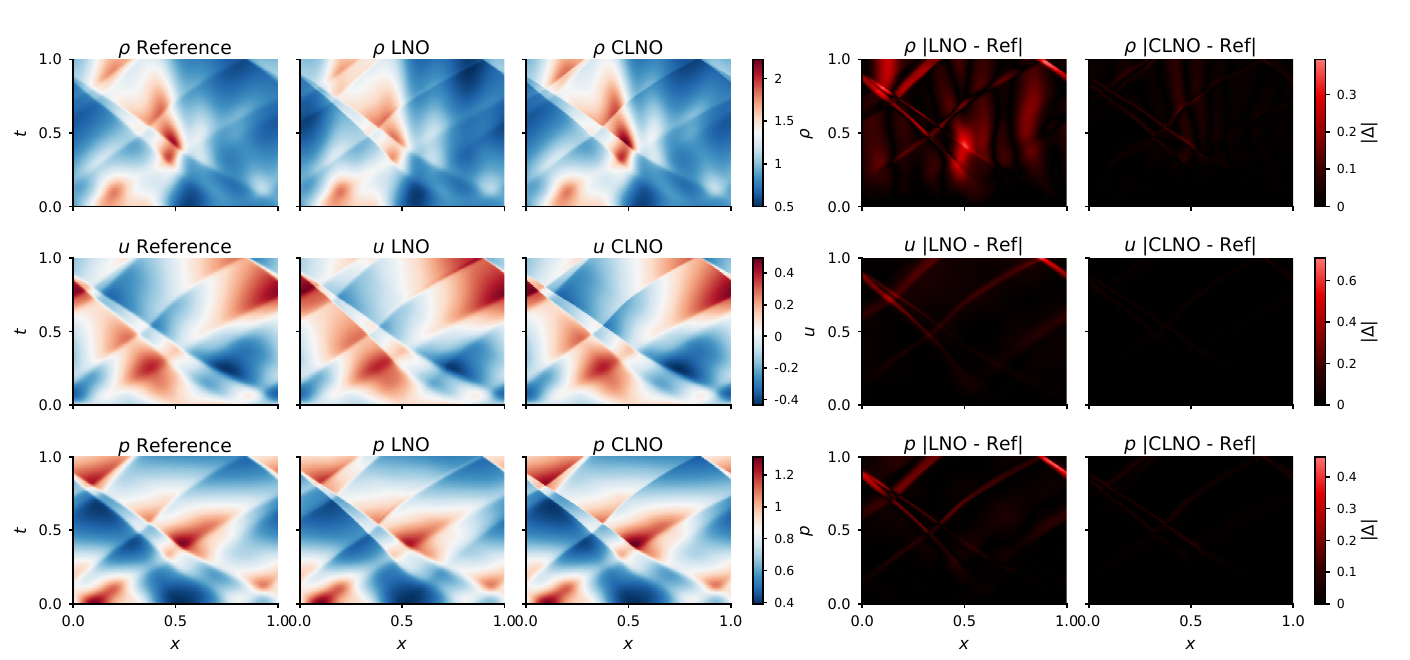}
\caption{Space--time evolution of density, velocity, and pressure
for the held-out periodic multiscale wave. Reference, LNO, CLNO,
and absolute-error fields are compared.}
\label{fig:one-d-fourier-fields}
\end{figure}

\begin{figure}[H]
\centering
\includegraphics[width=0.84\linewidth]{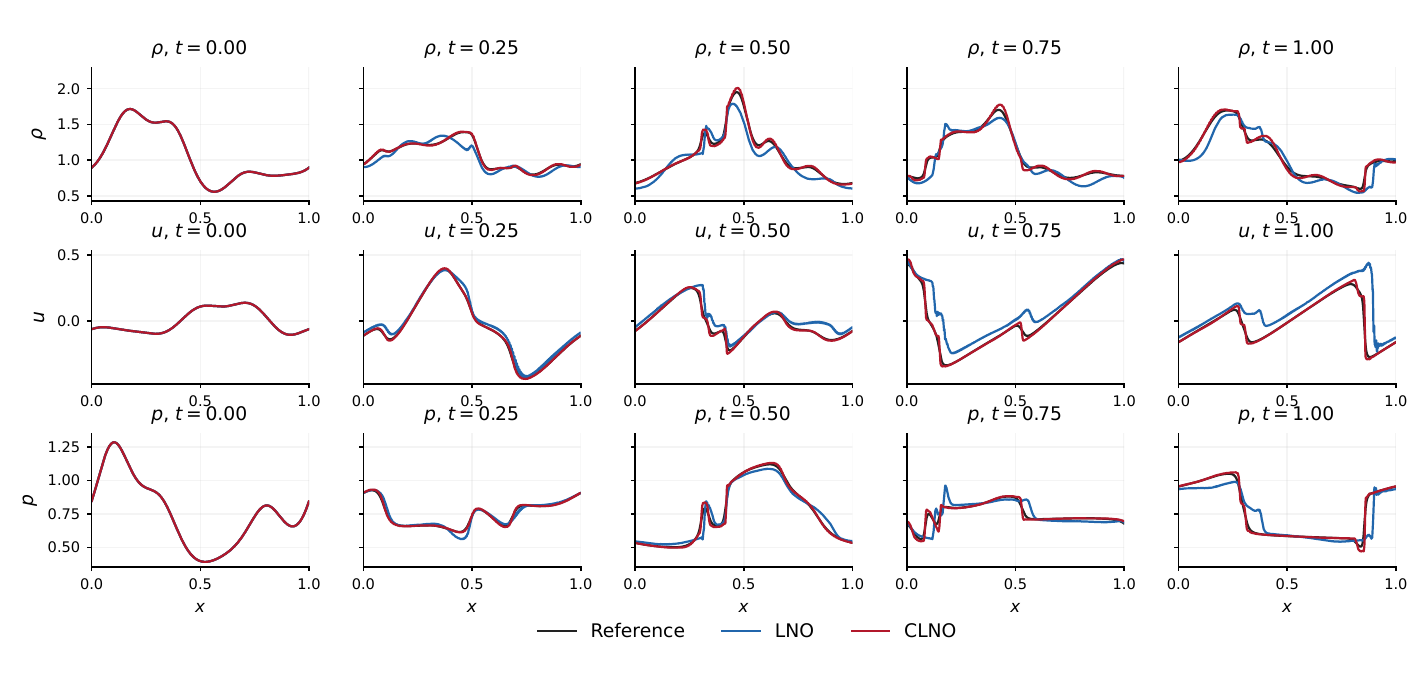}
\caption{Density, velocity, and pressure profiles of the periodic
multiscale wave at $t=0$, $0.25$, $0.5$, $0.75$, and $1$.}
\label{fig:one-d-fourier-snapshots}
\end{figure}

The periodic-wave and temporary-obstacle initial fields have initial
root-mean-square Mach number $\mathrm{Ma}_{\mathrm{rms}}(0)=0.089$;
the contact field has $0.181$. These describe the initial states and
are distinct from the prescribed Mach parameter in
Section~\ref{sec:parameter-generalization}. The reflecting solid
interval occupies $x/L\in[0.4375,0.5625]$ for $0.2\leq t<0.8$.

\begin{figure}[H]
\centering
\includegraphics[width=0.84\linewidth]{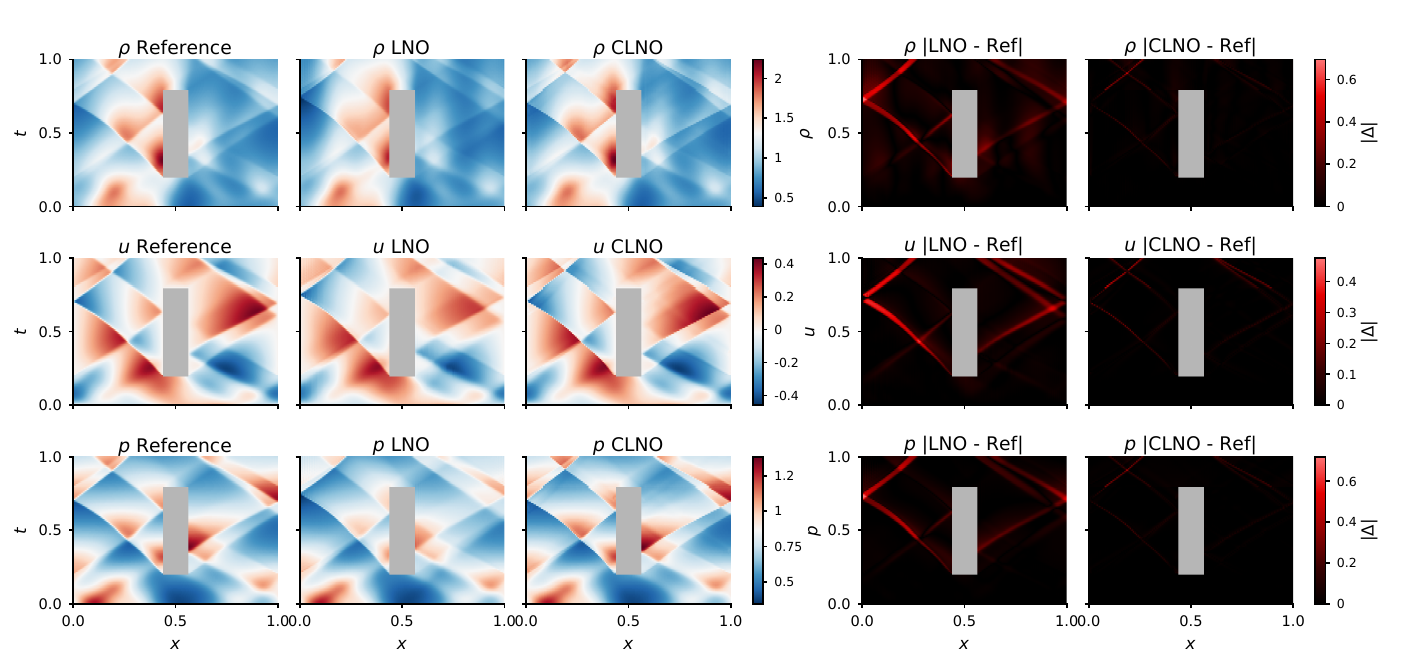}
\caption{Space--time density, velocity, and pressure for the
standard-grid temporary-obstacle test. The solid interval is
shaded gray during its active period.}
\label{fig:one-d-obstacle-fields}
\end{figure}

\begin{figure}[H]
\centering
\includegraphics[width=0.84\linewidth]{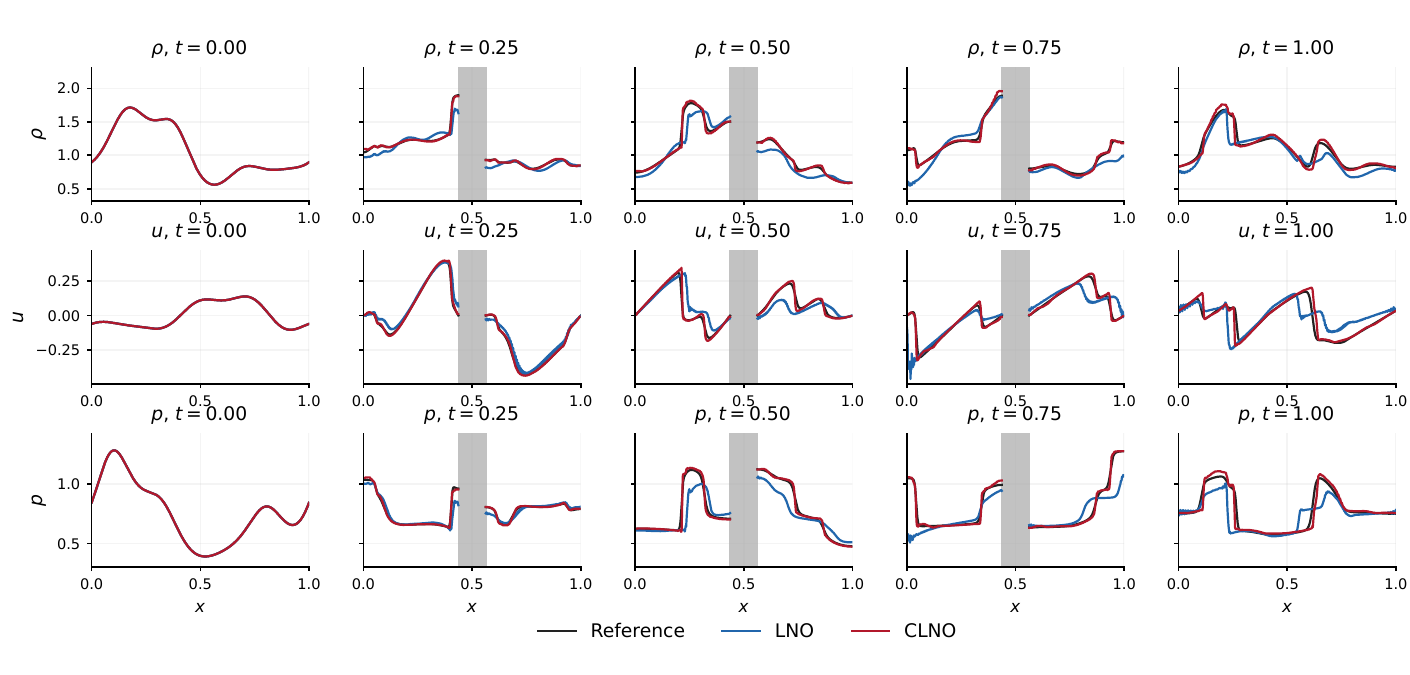}
\caption{Flow profiles with obstacle
appearance and removal at $t=0$,
$0.25$, $0.5$, $0.75$, and $1$.}
\label{fig:one-d-obstacle-snapshots}
\end{figure}

The nonperiodic contact test begins from
\begin{equation}
(\rho_0,u_0,p_0)=
\begin{cases}
(1.2,0.2,0.8),&x/L<0.53,\\
(0.6,0.2,0.8),&x/L\geq0.53,
\end{cases}
\label{app:eq:contact-initial}
\end{equation}
and uses constant-extrapolation outer boundaries. All reported
one-dimensional forecasts contain 100 future frames.

The additional pressure-discontinuity test uses the initial state
$(\rho_0,u_0,p_0)=(1,0,1)$ for $x<0.5$ and $(1,0,0.2)$ for
$x\geq0.5$ on $[0,1]$, with constant-extrapolation outer boundaries.
It uses the same $N=1024$, viscosity, reference discretization,
output interval, and 100-frame forecast protocol as the other
independent tests.

\begin{figure}[H]
\centering
\includegraphics[width=0.84\linewidth]{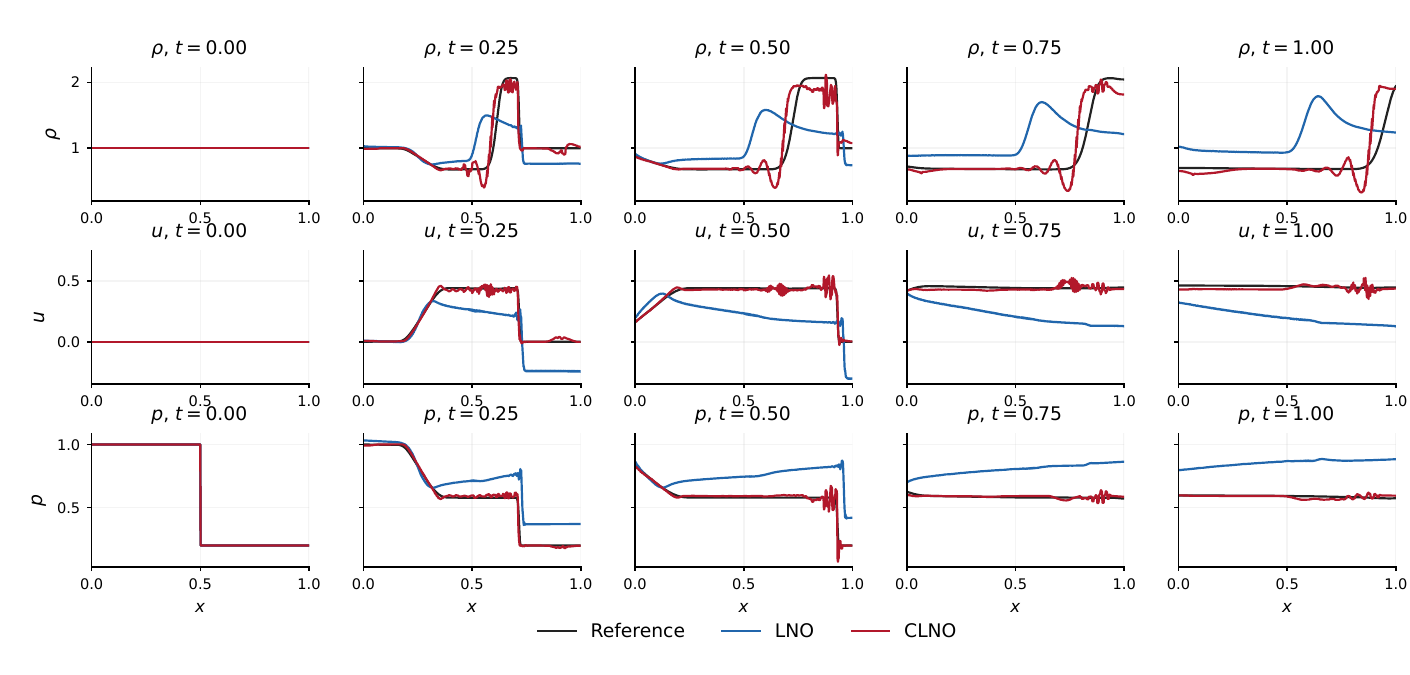}
\caption{Pressure-discontinuity profiles at $t=0$, $0.25$,
$0.5$, $0.75$, and $1$.}
\label{fig:one-d-pressure-snapshots}
\end{figure}

\begin{figure}[H]
\centering
\includegraphics[width=0.84\linewidth]{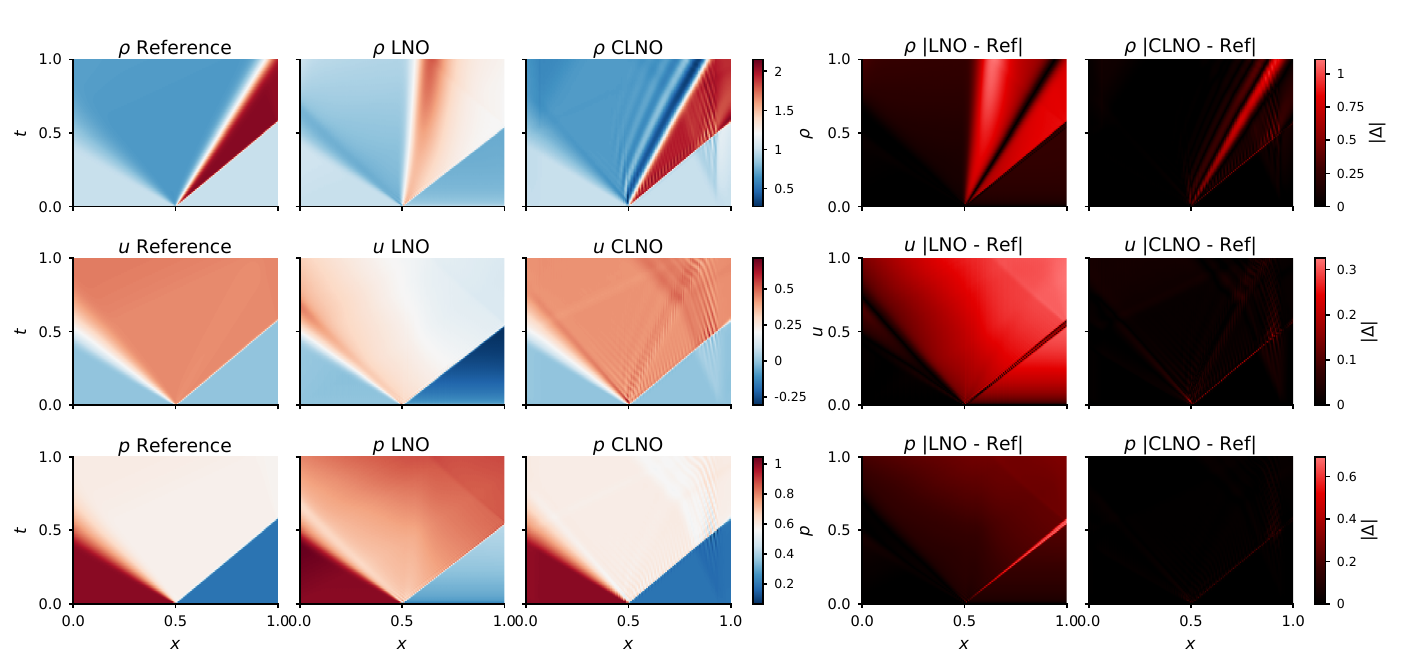}
\caption{Space--time density, velocity, and pressure for the
pressure-discontinuity Riemann problem. Columns show the reference,
LNO, CLNO, and their absolute errors.}
\label{fig:one-d-pressure-fields}
\end{figure}

\subsection{Fixed-parameter two-dimensional data}
\label{app:two-d-protocol}

The training domain is $[0,1]^2$, periodic in both directions,
with a $64\times64$ grid, $\gamma=1.4$, and viscosity $10^{-4}$
(nominal $\mathrm{Re}=10^4$). The reference solver combines
fifth-order WENO convective discretization and Lax--Friedrichs
flux splitting with fourth-order centered viscous terms and
third-order SSP Runge--Kutta integration. The 100 trajectories
contain 161 frames each on $t\in[0,8]$, with output interval
$\Delta t=0.05$. The training, validation, and held-out splits
contain 80, 10, and 10 trajectories. Stored primitive variables
$(\rho,u,v,p)$ have array shape $(100,4,64,64,161)$.

Training states contain three to five periodic vortices with
random positions, radii, and strengths, together with smooth
Fourier perturbations in density and pressure. The independent
high-wavenumber test uses modes through wave number eight; the
corresponding training perturbations use modes at most four.
The matched spatial-transfer test uses the same random seed and
initial-field generator in normalized coordinates $(x/L,y/L)$.
Its configurations are $(N,L)=(64,1),(96,1),(96,1.5)$, with
square domain $[0,L]^2$. The first two configurations use
$16\times16$ CLNO latent cells; the extended domain uses
$24\times24$ cells at the same latent-cell width. The equation
parameters, including viscosity, remain fixed. Each matched
forecast contains 20 frames at $\Delta t=0.05$.

The held-out and spatial-transfer cases start from a periodic
vortical field. Set \(\widetilde{\boldsymbol x}=(x/L,y/L)\) and let
\(\boldsymbol d_i\) be its wrapped displacement from the
\(i\)-th randomly sampled vortex center. Define
\(g_i=\exp[-|\boldsymbol d_i|^2/(2r_i^2)]\) and
\(\boldsymbol w=\sum_{i=1}^{N_v}(\omega_i g_i/r_i)(-d_{i,y},d_{i,x})\).
The initial state is
\begin{equation}
(\rho_0,\boldsymbol v_0,p_0)=\left(e^{0.054\mathcal W_\rho},\,\overline{\boldsymbol v}+0.3\frac{\boldsymbol w}{\langle|\boldsymbol w|^2\rangle^{1/2}},\,\rho_0^{\gamma_{\mathrm{gas}}}e^{0.015\mathcal W_{\mathrm{ent}}}/\gamma_{\mathrm{gas}}\right).
\label{eq:two-d-vortex-initial}
\end{equation}
Here \(N_v\in\{3,4,5\}\), \(r_i\in[0.07,0.16]\),
\(\omega_i\in[-1.25,-0.65]\cup[0.65,1.25]\), and
\(\overline{\boldsymbol v}\in[-0.12,0.12]^2\).
The ratio of specific heats is \(\gamma_{\mathrm{gas}}=1.4\).
The periodic Fourier fields are normalized to have a
root-mean-square (rms) value of one and take the form
\begin{equation}
\mathcal W_j(\widetilde{\boldsymbol x})=\frac{1}{\sigma_j}\sum_{0<|\boldsymbol k|\leq K}a_{j\boldsymbol k}|\boldsymbol k|^{-5/2}\cos(2\pi\boldsymbol k\cdot\widetilde{\boldsymbol x}+\phi_{j\boldsymbol k}),
\label{eq:two-d-fourier-initial}
\end{equation}
where \(a_{j\boldsymbol k}\) are standard-normal coefficients,
\(\phi_{j\boldsymbol k}\) are uniform phases, and
\(\sigma_j\) is the rms of the unnormalized sum.
The vortical cases use \(K=4\).
The held-out realization uses seed 20261020; the three matched
spatial configurations use seed 20261005 and the same normalized
coordinate field.

The independent high-wavenumber test uses seed 20261006 and
unit-rms Fourier fields \(\mathcal W_j\) with a maximum wavenumber of eight.
Writing \(\Phi=\mathcal W_\Phi\), \(\Psi=\mathcal W_\Psi\), and
\(a_{\mathrm{hf}}=0.3021276\), define
\(\boldsymbol q_{\mathrm{hf}}=(\mathcal D_x\Phi+\mathcal D_y\Psi,\mathcal D_y\Phi-\mathcal D_x\Psi)\).
Its primitive initial state is
\begin{equation}
(\rho_0,\boldsymbol v_0,p_0)=\left(e^{0.45a_{\mathrm{hf}}\mathcal W_\rho},\,\overline{\boldsymbol v}+\frac{a_{\mathrm{hf}}}{s_{\mathrm{hf}}}\boldsymbol q_{\mathrm{hf}},\,e^{0.35a_{\mathrm{hf}}\mathcal W_p}/\gamma_{\mathrm{gas}}\right),
\label{eq:two-d-highfreq-initial}
\end{equation}
where the mean velocity \(\overline{\boldsymbol v}\) is sampled from
\([-0.18,0.18]^2\),
\(\mathcal D_x,\mathcal D_y\) are the fourth-order
periodic difference operators, and \(s_{\mathrm{hf}}\) is the rms
magnitude of \(\boldsymbol q_{\mathrm{hf}}\). All these cases use
\(U(0,y,t)=U(L,y,t)\) and
\(U(x,0,t)=U(x,L,t)\), with periodic derivatives.

\begin{figure}[H]
\centering
\includegraphics[width=0.72\linewidth]{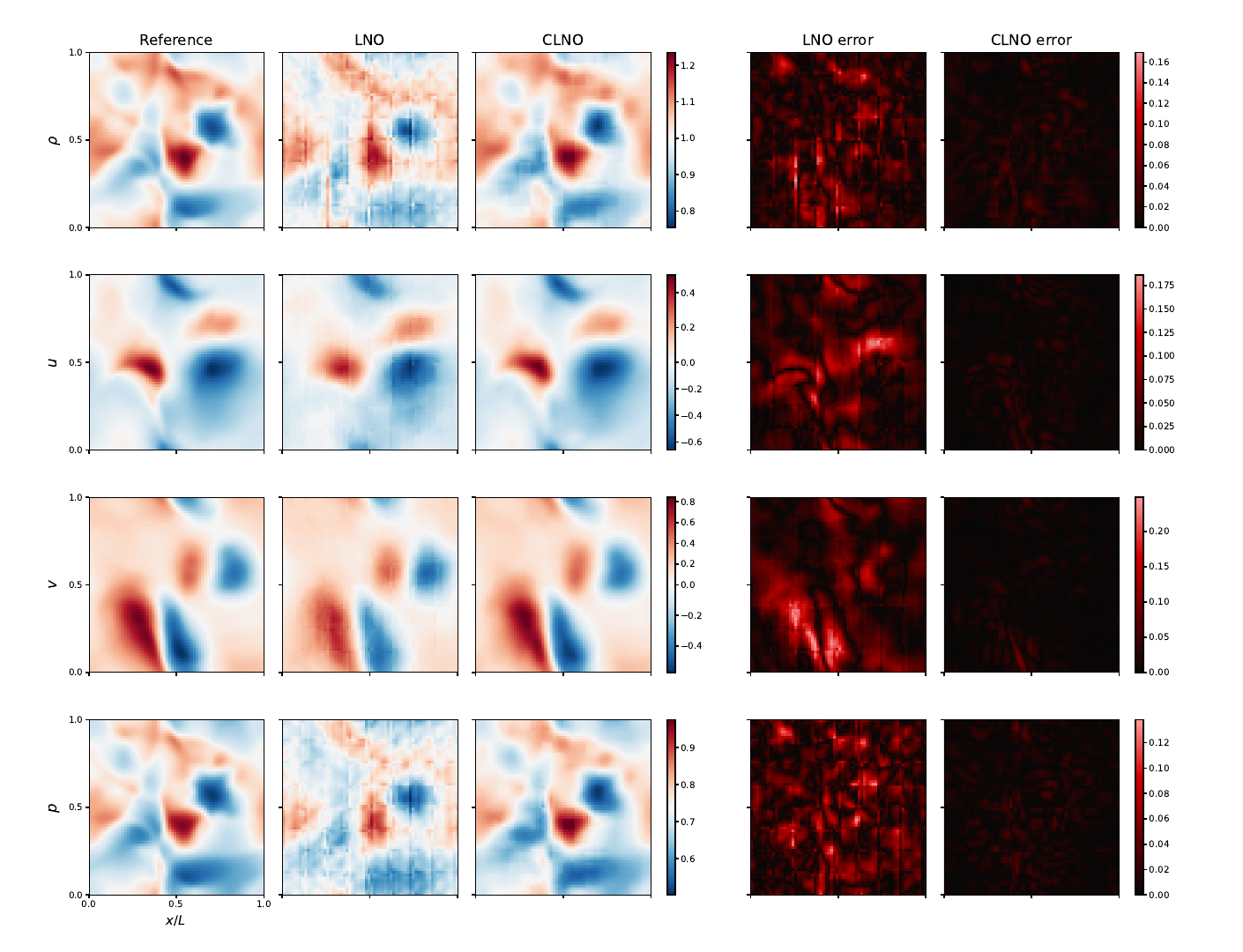}
\caption{Matched vortical flow on the standard $64\times64$ grid
over $[0,1]^2$ at $t=0.75$. Columns show the reference, LNO,
CLNO, and their absolute errors, respectively.}
\label{fig:two-d-64-standard}
\end{figure}

\subsection{Neural architecture configurations}
\label{app:architecture-config}

Table~\ref{tab:app-architecture} lists the architecture settings
used for the fixed-parameter comparisons. CLNO employs three
latent message blocks. Its direct interface-flux head applies a
$1$-point lift, a depthwise $3$-point convolution, a pointwise
mixing layer, and a final $3$-point convolution. The reconstruction
head has two $3$-point convolutions (effective stencil width five).
The one-dimensional decoder uses four nonconstant Legendre modes
per conservative variable. The two-dimensional decoder uses
24 nonconstant tensor-product modes per variable, with one
additional scalar amplitude gate for density and internal energy.
Both CLNO checkpoints use two internal substeps per stored frame.

\begin{table}[t]
\centering
\caption{Architecture settings of the reported fixed-parameter models.
$n,m,k$ are the local-resolution parameters of the LNO implementation.
FNO modes count retained Fourier modes per spatial direction;
DeepONet rank is the branch--trunk contraction dimension.}
\label{tab:app-architecture}
\small
\setlength{\tabcolsep}{4pt}
\renewcommand{\arraystretch}{1.08}
\begin{tabular}{@{}llll@{}}
\toprule
Benchmark & Model & Feature width & Spatial representation / depth\\
\midrule
1D NS & CLNO & 256 & 128 latent cells; 3 message blocks; 4 Legendre modes\\
1D NS & LNO & 20 & 4 blocks; $(n,m,k)=(12,6,2)$\\
1D NS & FNO & 40 & 4 spectral blocks; 16 Fourier modes\\
1D NS & DeepONet & 64 & branch pool 32; trunk rank 64\\
2D NS & CLNO & 256 & $16\times16$ latent cells; 3 blocks; 24 modes\\
2D NS & LNO & 40 & 4 blocks; $(n,m,k)=(16,8,2)$\\
2D NS & FNO & 40 & 4 spectral blocks; 12 modes per direction\\
2D NS & DeepONet & 64 & branch pool $8\times8$; trunk rank 64\\
\bottomrule
\end{tabular}
\end{table}

For the Reynolds--Mach study, the parameter-conditioned CLNO
retains the dimension-matched fixed-parameter widths, latent
geometry, and flux/reconstruction heads, and adds three FiLM
feature mappings at encoding, flux prediction, and reconstruction.
The one-dimensional FiLM CLNO checkpoint is evaluated with two
internal substeps per frame.

\subsection{Optimization and inference settings}
\label{app:training-inference}

Table~\ref{tab:app-training} reports the optimization settings
of the models used in the main tables. The one-dimensional LNO uses its separately trained
official baseline protocol; thus its optimizer and curriculum
differ from those of CLNO, FNO, and DeepONet. For the
$64\times64$ Navier--Stokes benchmark, CLNO and LNO share the
same 120-epoch curriculum and trajectory split.
Epoch counts give the minimum completed training length for each
row. Checkpoints were selected by validation error within each run.

\begin{table}[t]
\centering
\caption{Training settings for the reported checkpoints.
The 2D NS rows report a common minimum of 100 epochs;
the corresponding runs continued to 101 or 120 epochs.
Final horizon is the last curriculum horizon used for checkpoint
selection, rather than the number of frames in every evaluation.}
\label{tab:app-training}
\small
\setlength{\tabcolsep}{4pt}
\renewcommand{\arraystretch}{1.08}
\begin{tabular}{@{}llrrrr@{}}
\toprule
Benchmark & Model & Min. epochs & Batch & Initial LR & Final horizon\\
\midrule
1D NS & CLNO & 50 & 16 & $3\times10^{-4}$ & 10\\
1D NS & FNO / DeepONet & 50 & 16 & $3\times10^{-4}$ & 10\\
1D NS & LNO & 200 & 4 & $10^{-3}$ & 10\\
1D Re--Ma & CLNO (w/ FiLM) & 50 & 16 & $3\times10^{-4}$ & 10\\
1D Re--Ma & FNO / DeepONet & 50 & 16 & $3\times10^{-4}$ & 10\\
1D Re--Ma & LNO & 200 & 4 & $10^{-3}$ & 10\\
2D NS & CLNO / LNO & 100 & 2 & $3\times10^{-4}$ & 64\\
2D NS & FNO / DeepONet & 100 & 2 & $3\times10^{-4}$ & 64\\
2D NS & CLNO (w/o conservative flux) & 100 & 2 & $3\times10^{-4}$ & 64\\
2D NS & CLNO (w/o fine reconstruction) & 100 & 2 & $3\times10^{-4}$ & 64\\
\bottomrule
\end{tabular}
\end{table}

The one-dimensional CLNO curriculum is
$(1,2,4,8,10)$ frames for $(5,5,5,5,30)$ epochs,
respectively. Its optimizer is AdamW with weight decay $10^{-5}$,
gradient clipping at 1, and cosine learning-rate decay. The
two-dimensional fixed-parameter curriculum progresses through
$(1,2,4,8,12,16,24,32,40,48,64)$ frames, with ten epochs
per phase except 20 epochs at 64 frames. Both CLNO and LNO
use AdamW, weight decay $10^{-5}$, and cosine decay.
The w/o-conservative-flux ablation uses its epoch-101 checkpoint,
the first at the 64-frame horizon. The w/o-fine-reconstruction
ablation completed 120 epochs, with the epoch-117 checkpoint
selected by validation error.
For the fixed-parameter ablations, the flux ablation replaces
the shared interface-flux divergence with a learned per-cell
residual update, while the reconstruction ablation broadcasts
each updated latent-cell mean over its physical cells. Each
retains the remaining full-model components and uses the same
$64\times64$ training split and rollout curriculum.

Table~\ref{tab:app-inference} gives the numerical rollout
settings underlying the comparisons. A stored-frame step
$\Delta t$ is divided equally between CLNO's two internal
substeps. LNO, FNO, and DeepONet apply one learned map per
stored frame. The standard and spatial-transfer configurations
of each benchmark use the same trained weights.

\begin{table}[t]
\centering
\caption{Inference settings. The frame count is the number of
autoregressively predicted frames; all trajectories start at the
reference initial frame.}
\label{tab:app-inference}
\small
\setlength{\tabcolsep}{4pt}
\renewcommand{\arraystretch}{1.08}
\begin{tabular}{@{}lllll@{}}
\toprule
Benchmark & Reference interval & Forecast frames & CLNO substeps & Boundary\\
\midrule
1D fixed NS & $0.01$ & 100 & 2 & periodic / test-specific\\
2D fixed NS & $0.05$ & 20 or 160 & 2 & periodic\\
1D Re--Ma & $0.01$ & 64 & 2 & periodic\\
\bottomrule
\end{tabular}
\end{table}

\subsection{Parameter-dependent reference data}
\label{app:parameter-protocol}

The parameter study uses periodic finite-volume Rusanov
references with Newtonian stress and Fourier heat flux,
separate from the fixed-parameter WENO benchmarks.
Ten pairs combine $\mathrm{Re}\in\{5500,6500,7500,8500,9500\}$
with $\mathrm{Ma}\in\{0.35,0.45\}$ and reuse ten random
Fourier initial fields, giving 100 sequences. The pair
$(7500,0.45)$ is held out; the other nine provide 81 training
and nine validation sequences. The grid has 1024 cells on
$[0,1]$, without domain transfer.

For initial-field index \(i\), let \(\lambda_i\in[0.75,1.25]\) and
\(\mathcal H_q(x)=\frac{a_q}{2}\sum_{k=1}^{2}
[a_{qk}\sin(2\pi kx)+b_{qk}\cos(2\pi kx)]\), where the
coefficients are sampled independently from a standard normal
distribution and
\((a_\rho,a_u,a_p)=(0.08\lambda_i,0.20\lambda_i,
0.08\lambda_i)\). Every
parameter pair reuses the same ten primitive initial fields,
\begin{equation}
(\rho_0,u_0,p_0)=(\max\{1+\mathcal H_\rho,0.75\},\,\mathcal H_u,\,\max\{1+\mathcal H_p,0.75\}),\qquad U(0,t)=U(1,t).
\label{eq:re-ma-initial-main}
\end{equation}
The spatial derivatives also satisfy periodic matching at
\(x=0\) and \(x=1\).

All four methods use deterministic CUDA kernels with TF32 disabled in
Table~\ref{tab:re-ma-results}.

\clearpage
\addtocounter{page}{-1}
\end{document}